\documentclass[11pt]{article}
\usepackage{amsmath}
\usepackage{amssymb}
\usepackage{amsfonts}
\usepackage{cite}
\usepackage{amsmath,amsthm,amssymb}
\usepackage{mathrsfs}
\usepackage{color}
\allowdisplaybreaks[1]

\makeatletter

\newcommand{\Rmnum}[1]{\expandafter\@slowromancap\romannumeral #1@}
\makeatother

\theoremstyle{plain}
\newtheorem{theorem}{Theorem}[section]
\newtheorem{remark}[theorem]{Remark}
\newtheorem{lemma}[theorem]{Lemma}

\def\disp{\displaystyle}
\def\a{\alpha}

\def\t{\tau}

\def\l{\leq}

\def \B{\begin{equation}}
\def\bbb{\begin{array}{lll}}

\def\eee{\end{array}}
\def \E{\end{equation}}

\newif \ifLastSection \LastSectionfalse

\numberwithin{equation}{section}
\newcommand{\var}{\varepsilon}
\newcommand{\R}{\mathbb{R}}

\title{{Vacuum and singularity formation for compressible Euler equations with time-dependent damping
\thanks{This work is supported in part by the National Natural Science Foundation of China (Grant No. 11671237)}}}

\author{{Ying Sui $^a$\thanks{Email: suiying4320@163.com}}~~~~{Weiqiang Wang $^b$\thanks{Email: wangweiqiang@amss.ac.cn}}~~~~{Huimin Yu $^c$\thanks{Corresponding author, Email: hmyu@sdnu.edu.cn}}\\
\small\textit{$^{a,c}$ School of Mathematics and Statistics, Shandong Normal University, Jinan 250014, China.}\\
\small\textit{$^b$ Academy of Mathematics and Systems Science, Chinese Academy of Sciences, Beijing 100190, China.}\\}

\begin{document}
\date{}
\maketitle

\textbf{Abstract:}
In this paper, vacuum and singularity formation are considered for compressible Euler equations with time-dependent damping.
For $1<\gamma\l 3$, by constructing some new control functions ingeniously, we obtain the lower bounds estimates on density for arbitrary classical solutions.
Basing on these lower estimates, we succeed in proving the singular formation theorem for all $\lambda$, which was open in \cite{Sui-Yu} for some cases.
%In other words, our results not only include  the results in  \cite{Sui-Yu}, but also give the proofs to some unsolved cases in \cite{Sui-Yu}.
Moreover, the singularity formation of the compressible Euler equations when $\gamma=3$ is investigated, too. \\


\section{Introduction}
\

We consider the one dimensional compressible Euler equations with time-dependent damping in Lagrangian coordinates:
\begin{equation}\label{1.1}
	\left\{\begin{aligned}
		&\tau_{t}-u_{x}=0,\\
		&u_t+p_x=-\frac{\alpha}{(1+t)^{\lambda}}u,
	\end{aligned}
\right.
\end{equation}
with the initial data
\begin{equation}\label{1.2}
	\tau(x,0)=\tau_0(x),\;u(x,0)=u_0(x),
\end{equation}
where $\tau(x,t)$ and $u(x,t)$ are the specific volume and velocity of the flow respectively at location $x\in \R$ and time $t\in \R_{+}$. For simplicity, we assume the gas is deal polytropic and the gas pressure
\begin{equation}\label{1.3}
	p=K\tau^{-\gamma},~~\text{~for~}K>0\text{~and~}\gamma>1,
\end{equation}
where $\gamma$ is called  adiabatic index. Besides, $\alpha\geq 0, \ \lambda\in \R$ are two constants, and the term $-\frac{\alpha}{(1+t)^{\lambda}}u$ is the so-called damping effect on the fluid when $\alpha>0$.

For $\alpha=0$, the system \eqref{1.1} is the classical compressible Euler equation, which are the most basic model of hyperbolic conservation law equations and can be used to describe many physical phenomena. It is well known that, due to nonlinearity, classical solutions of hyperbolic conservation laws, may form gradient blowup in finite time, even when initial data is smooth and small. The study on breakdown of classical solutions for hyperbolic conservation laws has a long history \cite{Dafermos}. One can trace back to Stokes in \cite{Stokes} for a breakdown example for some scalar equation. In 1964, Lax \cite{Lax} studied the formation of singularities in finite time for $2\times 2$ general systems of strictly hyperbolic conservation laws when initial data is small and contains some compression, then followed by John, Li and Liu in \cite{John,Li,Liu} and etc., for more general systems of conservation laws. For $\alpha>0$, the system \eqref{1.1} can be used to describe the flow of fluids through porous media. However, the time-dependent damping makes some fantastic variety of the dynamic system. Therefore, the study on this problem is full of challenges and quite incomplete. Due to the practicability and strong relation with Euler equation, it has attracted extensive attention of many scholars in recent years, see \cite{2,3,4,5,6,7,8,9,10,11} for example.

%%Compressible Euler equations are the most basic model of hyperbolic conservation law equations, which can be used to describe many physical phenomena. It is well known that, due to nonlinearity, classical solutions of hyperbolic conservation laws, may form gradient blowup in finite time, even when initial data is smooth and small. The study on breakdown of classical solutions for hyperbolic conservation laws has a long history \cite{Dafermos}. One can trace back to Stokes in \cite{Stokes} for a breakdown example for some scalar equation. In 1964, Lax \cite{Lax} studied the formation of singularities in finite time for $2\times 2$ general systems of strictly hyperbolic conservation laws when initial data is small and contains some compression, then followed by John, Li and Liu in \cite{John,Li,Liu} and etc., for more general systems of conservation laws.

%%Due to the importance of compressible Euler equations in physics and the challenges be brought to mathematics, compressible Euler equations and its related problems have attracted extensive attention of many scholars, see \cite{1,2,3,4,5,6,7,8,9,10} for example.
 %%And compressible Euler system with time-dependent damping is an interesting topic and possesses rich wave phenomena. The study is also very technical and challenging but quite incomplete. Due to the technical reason, more attention has been paid in recent years.

In this paper, we study the vacuum and singularity formation of compressible Euler equations with time-dependent damping.
In particular, we extend the results and give the proofs of some unsolved cases in \cite{Sui-Yu}.
That is, we consider the cases of $\lambda\in \left(\min\{1,\frac{\alpha(\gamma-1)}{\gamma-3}\},\max\{1,\frac{\alpha(\gamma-1)}{\gamma-3}\}\right)$ when $ \gamma>3$ and $\lambda<\frac{\alpha(\gamma-1)}{\gamma-3}$ when  $1<\gamma<3$,  which are all open in \cite{Sui-Yu}.
 We also investigate the case of $\gamma=3$ and give the corresponding breakdown results.  Here, we still use the idea of constructing decoupled Riccati type equations
$$y'=a_0-a_2y^2
$$
for smooth solutions of compressible Euler equations. The difference between the two papers is that $a_0$ can change sign in this paper, while $a_0$ must be positive or negative in \cite{Sui-Yu}. Exactly speaking, for the unsolved case in \cite{Sui-Yu}, there exists a $t_0>0$ such that $a_0<0$ for $t>t_0$, while $a_0>0$ for $0 <t <t_0$. To overcome this difficulty, we firstly show the time-dependent lower bound on density for all $\lambda\in \R$ by constructing some new functions. The lower bound on density implies that $a_0$ is uniformly bounded for $0 <t <t_0$. Then using the uniform bound on $a_0$, we succeed in obtaining the singularity formation on the derivatives of velocity or specific volume if the initial compression reaches a certain level. Especially, our main results indicate that, for $1<\gamma\l3$, vacuum never occurs at any finite time provided the initial data away from vacuum, while the singularity must appear once the initial compression reaches a certain level, even for the over-damping case $\lambda< 0.$

 %For the unsolved case in \cite{Sui-Yu}, noticing the fact that there exists a $t_0>0$ such that $a_0<0$ for $t>t_0$, while $a_0$ is positive and bounded for $0 <t <t_0$, we can construct shock wave through a stronger initial compression. Then we succeed in proving the singularity formation theorem for all $\lambda$.
%In this paper, by constructing some new functions and obtain some sufficient conditions (only depending on the initial data at one point) to ensure the formation of singularity.

\section{The main Theorems}
\ \

For convenience, let's first introduce some variables and notations which have been used in \cite{Sui-Yu}. Denote
\begin{equation}\label{1.4}
\phi:=\int_{\tau}^{\infty}c(\tau)d\tau=\frac{2\sqrt{K\gamma}}{\gamma-1}\tau^{-\frac{\gamma-1}{2}}>0,
\end{equation}
where the nonlinear Lagrangian sound speed $c$ is defined as
\begin{equation}\label{1.4-1}
c:=\sqrt{-p_{\tau}}=\sqrt{K\gamma}\tau^{-\frac{\gamma+1}{2}}.
\end{equation}
It follows from \eqref{1.4}-\eqref{1.4-1} that
\begin{equation*}
	\tau=K_{\tau}\phi^{-\frac{2}{\gamma-1}},\ \ p=K_{p}\phi^{\frac{2\gamma}{\gamma-1}},  \ \ c=K_{c}\phi^{\frac{\gamma+1}{\gamma-1}},
\end{equation*}
where $K_{\tau}, K_{p}$ and $K_{c}$ are positive constants given by
\begin{equation*}
	K_{\tau}:=\left(\frac{2\sqrt{K\gamma}}{\gamma-1}\right)^{\frac{2}{\gamma-1}}, \ \
K_{p}:=KK^{-\gamma},\ \ \text{ and } \ \ K_{c}:=\sqrt{K\gamma}K_{\tau}^{-\frac{\gamma+1}{2}}.
\end{equation*}
%We also have
%$$
%K_{p}=\frac{\gamma-1}{2\gamma}K_{c}\text{ and } %K_{\tau}K_{c}=\frac{\gamma-1}{2}.
%$$

A direct calculation shows that the system \eqref{1.1} has two characteristic speeds
$$
\lambda_1=-\lambda_2=c.
$$
The forward and backward characteristic lines are described by
$$
\frac{dx_{+}(t)}{dt}=c\quad \text{ and }\quad  \frac{dx_{-}(t)}{dt}=-c,
$$
and we denote the corresponding directional derivatives along them by
$$
\prime=\frac{\partial}{\partial t}+c\frac{\partial}{\partial x}\quad \text{ and }\quad
\backprime=\frac{\partial}{\partial t}-c\frac{\partial }{\partial x}
$$
respectively. Furthermore, we denote the Riemann invariants by
\begin{equation}\label{1.5-1}
w:=u+\phi\quad \text{ and }\quad z:=u-\phi.
\end{equation}
By a direct calculation, we have
\begin{equation}\label{1.5}
w^\prime =-\frac{\alpha}{2(1+t)^{\lambda}}(z+w),\qquad
z^\backprime =-\frac{\alpha}{2(1+t)^{\lambda}}(z+w).
\end{equation}
%\begin{equation}\label{1.6}
%z^\backprime :=-\frac{\alpha}{2(1+t)^{\lambda}}(z+w)
%\end{equation}

In \cite{Sui-Yu}, the authors have proved the following theorems.

\begin{theorem}[{\cite[Theorem 3.2]{Sui-Yu}}]\label{thm1.2}
	Suppose the initial data $(\tau_0,u_0)(x)\in C^1(\R)$ and there exists a positive constant $C_0$ such that
	\begin{equation}\label{1.10}
	\|u_0(x)\|_{C^1}\l C_0,~ \|\tau_{0}(x)\|_{C^1}\l C_0, ~ C_0^{-1}< \tau_{0}(x).
	\end{equation}
	Let $1<\gamma<3$, $\lambda\geq \frac{\alpha(\gamma-1)}{\gamma-3}$ and $(\tau,u)$ is a $C^1$ solution of \eqref{1.1}-\eqref{1.2}. Assume there exists one point $x_0$ such that
	\begin{equation}\label{1.5-2}
	u_{x}(x_0,0)+\phi_x(x_0,0)<-\frac{\alpha(\gamma-1)}{K_{c}(3-\gamma)}(\phi(x_0,0))^{\frac{-2}{\gamma-1}}
	\end{equation}
	or
	\begin{equation}\label{1.5-3}
	u_{x}(x_0,0)-\phi_x(x_0,0)<-\frac{\alpha(\gamma-1)}{K_{c}(3-\gamma)}(\phi(x_0,0))^{\frac{-2}{\gamma-1}}.
	\end{equation}
	Then $\|(\tau_x, u_x)\|_{L^\infty}$ must blow up in finite time.
\end{theorem}

\begin{theorem}[{\cite[Theorem 3.1]{Sui-Yu}}]\label{thm1.1}
Suppose the initial data $(\tau_0,u_0)(x)\in C^1(\R)$ and there exists a positive constant $C_0$ such that
	\begin{equation}\label{1.7}
	|u_0(x)|\leq C_0,\quad C_0^{-1}<\tau_0(x).
	\end{equation}
Let $\gamma>3, \lambda<\min\{1,\frac{\alpha(\gamma-1)}{\gamma-3}\}$ or $\lambda>\max\{1,\frac{\alpha(\gamma-1)}{\gamma-3}\}$, and $(\tau,u)$ is a $C^1$ solution of \eqref{1.1}-\eqref{1.2}. Assume there exists one point $x_0$ such that
\begin{equation}\label{1.5-5}
u_{x}(x_0,0)+\phi_{x}(x_0,0)<K_{1}(\phi(x_0,0))^{-\frac{2}{\gamma-1}}-K_2(\phi(x_0,0))^{-\frac{\gamma+1}{2(\gamma-1)}}
\end{equation}
or
\begin{equation}\label{1.5-6}
u_{x}(x_0,0)-\phi_{x}(x_0,0)<K_{1}(\phi(x_0,0))^{-\frac{2}{\gamma-1}}-K_2(\phi(x_0,0))^{-\frac{\gamma+1}{2(\gamma-1)}},
\end{equation}
where $K_{1},K_{2}$ are positive constants which depend only on $C_0,\gamma,K,\lambda$ and $\alpha$. Then $\|(\tau_x, u_x)\|_{L^\infty}$ must blow up in finite time.
\end{theorem}

For the rest cases, that is,
\begin{equation}\label{a2}
	\lambda<\frac{\alpha(\gamma-1)}{\gamma-3} \text{ when } 1<\gamma<3,
\end{equation}
and
\begin{equation}\label{a1}
	\lambda\in \left(\min\{1,\frac{\alpha(\gamma-1)}{\gamma-3}\},\max\{1,\frac{\alpha(\gamma-1)}{\gamma-3}\}\right) \text{ when }\gamma>3,
\end{equation}
the authors assert that this is just a technical problem and the Riemann invariants may blow up for some certain initial data, too. Indeed, we consider the rest cases in this paper and have obtained the following theorems.

\begin{theorem}[Singularity formation for $1<\gamma<3$ and $\lambda\neq 1$]\label{thm1.4}
Suppose the initial data $(\tau_0,u_0)(x)\in C^1(\R)$ satisfying \eqref{1.10} and $(\tau,u)$ is a $C^1$ solution of \eqref{1.1}-\eqref{1.2}.
%there exists a positive constant $C_0>0$ such that
%\begin{equation}\label{1.10}
%|u_0(x)|_{C^1}\l C_0,~  |\tau_{0}(x)|_{C^1}\l C_0, ~ C_0^{-1}< \tau_{0}(x).
%\end{equation}
Then $\|u\|_{L^{\infty}}$ will be uniformly bounded for all time $t>0$ and the density $\rho=\tau^{-1}$ has a time-dependent lower bound,
\begin{equation}\label{1.20}
	\begin{aligned}
		&\rho(x,t)=\tau^{-1}(x,t)\geq {\left\{\begin{matrix} C\left(1+ t\right)^{\frac{4}{\gamma-3}}\quad\quad\quad\quad\quad&\quad \text{for }\lambda\geq0,
				\\C\left(1+ e^{\frac{-\lambda(1+t)^{1-\lambda}}{1-\lambda}}t\right)^{\frac{4}{\gamma-3}}&\quad  \text{for }\lambda<0,
			\end{matrix}
			\right.}
	\end{aligned}
\end{equation}
where $C$ is a positive constants depending only on $C_0, \gamma, K, \lambda$ and $\alpha$. \eqref{1.20} indicates that the vacuum will never occur in finite time. Moreover, assume there exists one point $x_0$ such that
	\begin{equation}\label{1.12}
		u_{x}(x_0,0)+\phi_{x}(x_0,0)<K_{3}(\phi(x_0,0))^{-\frac{2}{\gamma-1}}-K_4(\phi(x_0,0))^{-\frac{\gamma+1}{2(\gamma-1)}}
	\end{equation}
	or
	\begin{equation}\label{1.13}
		u_{x}(x_0,0)-\phi_{x}(x_0,0)<K_{3}(\phi(x_0,0))^{-\frac{2}{\gamma-1}}-K_4(\phi(x_0,0))^{-\frac{\gamma+1}{2(\gamma-1)}},
	\end{equation}
	where $K_3,K_4$ are positive constants which depend only on $C_0,\gamma,K,\lambda$ and $\alpha$. Then $\|(\tau_x, u_x)\|_{L^\infty}$ must blow up in finite time.
\end{theorem}

\begin{theorem}[Singularity formation for $\gamma>3$ and $\lambda\neq 1$]\label{thm1.3}
	Suppose the initial data $(\tau_0,u_0)(x)\in C^1(\R)$ satisfying \eqref{1.7} and $(\tau,u)$ be a $C^1$ solution of \eqref{1.1}-\eqref{1.2}.	
%and there exists a positive constant $C_0$ such that
%\begin{equation}\label{1.7}
%|u_0(x)|\leq C_0,\quad C_0^{-1}<\tau_0(x).
%\end{equation}
	Assume there exists one point $x_0$ such that
\begin{equation}\label{1.8}
	u_{x}(x_0,0)+\phi_{x}(x_0,0)<K_{5}(\phi(x_0,0))^{-\frac{2}{\gamma-1}}-K_6(\phi(x_0,0))^{-\frac{\gamma+1}{2(\gamma-1)}}
\end{equation}
	or
\begin{equation}\label{1.9}
	u_{x}(x_0,0)-\phi_{x}(x_0,0)<K_{5}(\phi(x_0,0))^{-\frac{2}{\gamma-1}}-K_6(\phi(x_0,0))^{-\frac{\gamma+1}{2(\gamma-1)}},
\end{equation}
where $K_5,K_6$ are two positive constants which depend only on $C_0,\gamma,K,\lambda$ and $\alpha$. Then $\|(\tau_x, u_x)\|_{L^\infty}$ must blow up in finite time.
\end{theorem}

\begin{remark}
	\eqref{1.20} holds for all $\lambda\in \R$, and the condition $\lambda\neq 1$ in Theorem \ref{thm1.4} is only used in proving the blow-up of $\|(\tau_{x},u_{x})\|_{L^{\infty}}$. Singularity formation for $\lambda=1$ with all $\gamma>1$ can be found in \cite[Theorems 4.1 and 4.2]{Sui-Yu}, also see the Remark \ref{rem4.1} below.
\end{remark}

%%%%%%%%%%%%%%%%%%%%%%%%%%%%%%%%%%%%%%%%%%%%%%%%%%%%%%%%%%%%%%%%%%%%%%%%%%%%%%%%%%%%%%%%%%%%%%%%%%%%%%%%%%%%%%%%%%%%%%%%%%%%%%%%%%%%%%%%%
%%%%%%%%%%%%%%%%%%%%%%%%%%%%%%%%%%%%%%%%%%%%%%%%%%%%%%%%%%%%%%%%%%%%%%%%%%%%%%%%%%%%%%%%%%%%%%%%%%%%%%%%%%%%%%%%%%%%%%%%%%%%%%%%%%%%%%%%%%%%%%

Theorem \ref{thm1.4} and Theorem \ref{thm1.3} remove the restrictions on $\lambda$ in \cite{Sui-Yu} and include the cases $(\ref{a2})$ and $(\ref{a1})$. In addition, the singularity formation of the compressible Euler equations when $\gamma=3$ is given in the following theorem.

\begin{theorem}[Singularity formation for $\gamma=3$]\label{thm9.1}
Suppose the initial data $(\tau_0,u_0)(x)\in C^1(\R)$ satisfying \eqref{1.10} and $(\tau,u)(x,t)$ is a $C^1$ solution of $\eqref{1.1}-\eqref{1.2}$ with $\gamma=3$.
%and there exists a constant $C_0>0$ such that
%\begin{equation}
    %|u_0(x)|_{C^1}\l C_0,~  |\tau_{0}(x)|_{C^1}\l C_0, ~ C_0^{-1}< \tau_{0}(x).
%\end{equation}
Then $\|u\|_{L^{\infty}}$ will be uniformly bounded for all time $t>0$ and the density $\rho=\tau^{-1}$ has a time-dependent lower bound,
\begin{equation}\label{1.200}
	\begin{aligned}
	&\rho(x,t)=\tau^{-1}(x,t)\geq {\left\{\begin{matrix} C_{0}^{-1}e^{-\hat{M}t}\quad\quad\quad\quad&\quad\quad \text{for }\lambda\geq0,
	\\C_{0}^{-1}e^{-\hat{M}te^{\frac{-\lambda(1+t)^{1-\lambda}}{1-\lambda}}}&\quad\quad  \text{for }\lambda<0,
	\end{matrix}
	\right.}
	\end{aligned}
\end{equation}
where $\hat{M}$ is a positive constant depending only on $C_0, K,\gamma,\lambda$ and $\alpha$. \eqref{1.200} indicates that the vacuum will never occur in finite time.
Moreover, for any given $t_0>0$, there exists a non-decreasing (w.r.t. $t_0$) function $\tilde{M}(t_0)$, depending on $C_0, K,\gamma,\lambda$ and $\alpha$, such that if there exists one point $x_0$ satisfying
\begin{equation}\label{9.14}
	\begin{aligned}
		u_{x}(x_0,0)+\phi_{x}(x_0,0)<\frac{\alpha}{2K_{c}}\phi^{-1}(x_0,0)\ln \phi(x_0,0)-\max\{\frac{2}{K_{c}t_0}, \tilde{M}(t_0)\},
	\end{aligned}
\end{equation}
or
\begin{equation}\label{9.15}
	\begin{aligned}
		u_{x}(x_0,0)-\phi_{x}(x_0,0)<\frac{\alpha}{2K_{c}}\phi^{-1}(x_0,0)\ln \phi(x_0,0)-\max\{\frac{2}{K_{c}t_0}, \tilde{M}(t_0)\}.
	\end{aligned}
\end{equation}
Then $\|u_x,\tau_x\|_{L^\infty}$ must blow up before $t_0$.
\end{theorem}

\begin{remark}
	We mention that we only need $\eqref{1.10}$ to obtain the time-dependent lower bound on density \eqref{1.20} for $1<\gamma<3$, and \eqref{1.200} for $\gamma=3$. Therefore, \eqref{1.20} and \eqref{1.200} may be viewed as a generalization of \cite[Lemma 3.3]{Sui-Yu} and \cite[Lemma 12]{11}, in which the authors impose some restrictions on $\lambda$ and $\alpha$.
\end{remark}

\begin{remark}
	For $1<\gamma\leq3$ and $\lambda=0$, we can get the lower bound of density independent of $t$.
In other words, using  the proof of Theorem \ref{thm1.4} and Theorem \ref{thm9.1}, we can respectively obtain the fixed lower bound of density for $1<\gamma\leq3$ and $\gamma=3$ by simple calculation.
\end{remark}

\begin{remark}
      The time-dependent lower bound on density plays an important role in obtaining blow-up of $\|(\tau_{x},u_{x})\|_{L^{\infty}}$ in Theorems \ref{thm1.4} and \ref{thm9.1}. For $\gamma>1$, under the assumption \eqref{1.10}, an optimal time-dependent lower bound on density is obtained in \cite{1,Chen1} for compressible isentropic Euler equations and in \cite{Chen2} for full compressible Euler equations. However, for $\gamma\geq 3$, as stated in \cite[Remark 3]{11}, it is open whether vacuum occurs in finite time or not for damped Euler equations. Theorem \ref{thm9.1} give a negative answer to this open question for $\gamma=3$, i.e., the vacuum will never occur in finite time.
\end{remark}

For latter use, we recall some useful lemmas, whose proof can be found in \cite{Sui-Yu}.

\begin{lemma}[{\cite[Theorem 2.1]{Sui-Yu}}]\label{lem1.1}
Let $(\tau_0,u_0)(x)\in C^{1}(\R)$ satisfying \eqref{1.7} and $(\tau,u)$ be the $C^1$ solution to system (\ref{1.1})-(\ref{1.2}). Then in the domain where the solution exists, $(\tau,u)(x,t)$ is uniform bounded in the following form
\begin{equation}\label{6191}
\begin{aligned}
|u(x,t)|\leq \tilde{C}_{0},\quad \tau(x,t)\geq \tilde{C}_{0}^{-1},
 \end{aligned}
\end{equation}
where $\tilde{C}_0$ can be any constant bigger than $\max\{C_0+C_{0}^{\frac{\gamma-1}{2}},(C_0+C_0^{\frac{\gamma-1}{2}})^{\frac{2}{\gamma-1}}\}$.
\end{lemma}

\begin{lemma}[{\cite[Lemma 2.2]{Sui-Yu}}]\label{L6.1}
	The $C^1$ solutions $(\tau,u)$ of (\ref{1.1})-(\ref{1.2}) satisfy
	\begin{equation*}
		\begin{aligned}
			{A}^\prime= -\frac{\a}{2(1+t)^\lambda}({A}+{B})+K_d({A}{B}-{A}^2)
		\end{aligned}
	\end{equation*}
	and
	\begin{equation*}
		\begin{aligned}
			{B}^\backprime=-\frac{\a}{2(1+t)^\lambda}({A}+{B})+K_d({A}{B}-{B}^2),
		\end{aligned}
	\end{equation*}
	where
	$$
	{A}={w}_x,\quad {B}={z}_x,\quad K_d =K_c\frac{\gamma+1}{2(\gamma-1)}\phi^{\frac{2}{\gamma-1}}.
	$$
\end{lemma}

\begin{lemma}[{\cite[Lemma 3.1]{Sui-Yu}}]\label{lem1.2}
For the $C^1$ solutions $(\tau,u)$ of (\ref{1.1})-(\ref{1.2}), denoting
\begin{equation}\label{1.14-1}
y:=\left(\phi^{\frac{\gamma+1}{2(\gamma-1)}}w_{x}-\frac{\alpha(\gamma-1)}{K_{c}(\gamma-3)(1+t)^{\lambda}}\phi^{\frac{\gamma-3}{2(\gamma-1)}}\right)e^{\frac{\alpha(3\gamma-1)}{2(\gamma-3)(1-\lambda)}(1+t)^{-\lambda}},
\end{equation}
and
\begin{equation}\label{1.14-2}
q:=\left(\phi^{\frac{\gamma+1}{2(\gamma-1)}}z_{x}-\frac{\alpha(\gamma-1)}{K_{c}(\gamma-3)(1+t)^{\lambda}}\phi^{\frac{\gamma-3}{2(\gamma-1)}}\right)e^{\frac{\alpha(3\gamma-1)}{2(\gamma-3)(1-\lambda)}(1+t)^{-\lambda}},
\end{equation}
then $y$ and $q$ satisfy
\begin{equation}\label{1.14}
	  y^\prime=a_0-a_2y^2,\quad   q^\backprime=a_0-a_2q^2,
\end{equation}
where
\begin{equation}\label{1.16}	a_2=\frac{K_{c}(\gamma+1)}{2(\gamma-1)}\phi^{-\frac{\gamma-3}{2(\gamma-1)}}e^{-\frac{\alpha(3\gamma-1)}{2(\gamma-3)(1-\lambda)}(1+t)^{1-\lambda}}>0,
\end{equation}

\begin{equation}\label{1.15}
	a_0=\frac{\lambda\alpha(\gamma-1)(\gamma-3)(1+t)^{\lambda-1}-\alpha^2(\gamma-1)^2}{K_{c}(\gamma-3)^2(1+t)^{2\lambda}}\phi^{\frac{\gamma-3}{2(\gamma-1)}}e^{\frac{\alpha(3\gamma-1)}{2(\gamma-3)(1-\lambda)}(1+t)^{1-\lambda}}.
\end{equation}
\end{lemma}

%\begin{lemma}[{\cite[Lemma 3.2]{Sui-Yu}}]\label{L1.11}
% Suppose $1<\gamma<3,$ $ \ld\g {{\a(\gamma-1)}\over{\gamma-3}}$. Then for any $C^1$ solutions of  (\ref{1.1})-(\ref{1.2}), we have a priori bounds
%\begin{equation}\label{1.18}
%\begin{aligned}
%y(x, t)\leq \max\{1, \sup_x(y(x,0))\}:=Y,
%\end{aligned}
%\end{equation}
%and
%\begin{equation}\label{1.19}
%\begin{aligned}
%q(x,t)\leq  \max \{1, \sup_x(q(x,0))\}:=Q.
%\end{aligned}
%\end{equation}
%\end{lemma}

\ \
\section{Proof of the main theorems}
\subsection{Proof of Theorem \ref{thm1.4}}
\ \

We divide the proof into five steps. In steps 1-4, we first give the proof of the lower bound on density, that is, the proof of \eqref{1.20}. In step 5, we prove the blow-up of $\|(\tau_{x},u_{x})\|_{L^{\infty}}$. \\

\textbf{Step 1.} It follows from \eqref{1.5} and the fact $c_{w}=-c_{z}$ that
\begin{equation}\label{4.6}
	\begin{aligned}
		(w_x)'&=(w')_x-c_xw_x\\
		&=-\frac{\a}{2(1+t)^\lambda}({z}+{w})_x-(c_ww_x^2+c_zz_xw_x)\\
		&=-\frac{\a}{2(1+t)^\lambda}({z}+{w})_x-c_ww_x^2+c_wz_xw_x,
	\end{aligned}
\end{equation}
and
\begin{equation}\label{4.7}
	\begin{aligned}
		(z_x)^\backprime&=(z^\backprime)_x+c_xz_x\\
		&=-\frac{\a}{2(1+t)^\lambda}({z}+{w})_x+(c_ww_xz_x+c_zz_x^2)\\
		&=-\frac{\a}{2(1+t)^\lambda}({z}+{w})_x-c_wz_x^2+c_wz_xw_x.
	\end{aligned}
\end{equation}
Noticing from \eqref{1.4}-\eqref{1.4-1}, it holds
$$
	\begin{aligned}
		\phi_\t=-c,
	\end{aligned}
$$
which together with \eqref{1.1}, yields
$$
\begin{aligned}
	\phi^\prime&=\phi_t+c\phi_x
	=-c\t _t+c\phi_x=-cu_x+c\phi_x=-c(u-\phi)_x=-cz_x,\\
	\phi^\backprime&=\phi_t-c\phi_x
	=-c\t _t-c\phi_x=-cu_x-c\phi_x=-c(u+\phi)_x=-cw_x.
\end{aligned}
$$
Thus, we have
\begin{equation}\label{4.8-1}
%\begin{aligned}
	z_x=-\frac{1}{c}\phi^\prime=-\frac{1}{c}\phi_\tau\tau^\prime=\tau^\prime,\qquad
w_{x}=-\frac{1}{c}\phi^\backprime=-\frac{1}{c}\phi_\tau\tau^\backprime=\tau^\backprime.
%\end{aligned}
\end{equation}
Moreover, using \eqref{4.8-1}, a direct calculation shows that
\begin{equation}\label{4.8-2}
\begin{aligned}
	\frac{1}{2}c^{-\frac{1}{2}}c'w_x&=\frac{1}{2}c^{-\frac{1}{2}}(c_ww_\t\t'+c_zz_\t\t')w_x
	=c^{-\frac{1}{2}}c_ww_\t\t'w_x=-c^{\frac{1}{2}}c_wz_xw_x,
\end{aligned}
\end{equation}
where we have used $c_w=-c_z$, $w_\t=-z_\t$ and $w_{\tau}=\phi_{\tau}=-c$.

Hence it follows form \eqref{4.6}-\eqref{4.7}, \eqref{4.8-2} that
\begin{equation}\label{4.9}
	\begin{aligned}
		(c^{\frac{1}{2}}w_x)'&= c^{\frac{1}{2}}(w_x)'+\frac{1}{2}c^{-\frac{1}{2}}c'w_x\\
		&=-\frac{\a c^{\frac{1}{2}}}{2(1+t)^\lambda}({z}+{w})_x-c^{\frac{1}{2}}c_ww_x^2+c^{\frac{1}{2}}c_wz_xw_x+\frac{1}{2}c^{-\frac{1}{2}}c'w_x.\\
		&=-\frac{\a c^{\frac{1}{2}}}{2(1+t)^\lambda}({z}+{w})_x-c^{\frac{1}{2}}c_ww_x^2,
	\end{aligned}
\end{equation}
and
\begin{equation}\label{4.9-1}
\begin{aligned}
	(c^{\frac{1}{2}}z_x)^\backprime&= c^{\frac{1}{2}}(z_x)^\backprime+\frac{1}{2}c^{-\frac{1}{2}}c^\backprime z_x\\
	&=-\frac{\a c^{\frac{1}{2}}}{2(1+t)^\lambda}({z}+{w})_x-c^{\frac{1}{2}}c_wz_x^2+c^{\frac{1}{2}}c_wz_xw_x
	+\frac{1}{2}c^{-\frac{1}{2}}c^\backprime z_x\\
	&=-\frac{\a c^{\frac{1}{2}}}{2(1+t)^\lambda}({z}+{w})_x-c^{\frac{1}{2}}c_wz_x^2.
\end{aligned}
\end{equation}
%Besides
%by (\ref{4.8}), we have $w_\t=\phi_\t=-c$, then $\t_w=\phi_\t^{-1}=-\frac{1}{c}$.
%Then (\ref{4.9}) becomes
%\begin{equation}\label{4.10}
	%(c^{\frac{1}{2}}w_x)'
	%=-\frac{\a %c^{\frac{1}{2}}}{2(1+t)^\lambda}({z}+{w})_x-c^{\frac{1}{2}}c_ww_x^2.
%\end{equation}
%Similarly, from (\ref{4.7}), we have
\\

\textbf{Step 2.} For $\lambda\geq 0$ and $\lambda\neq 1$, we introduce two new variables
\begin{equation}\label{b.1}
	G(t,x)=c^{1\over 2}w_x+{\alpha\over {(1+t)^\lambda}}h(\tau),
\end{equation}
\begin{equation}\label{b.2}
	H(t,x)=c^{1\over 2}z_x+{\alpha\over {(1+t)^\lambda}}h(\tau),
\end{equation}
where
\begin{equation}\label{b.02}
	h(\t) =\int^{\t}_{0}c(s)^{\frac{1}{2}}ds=\frac{4(K\gamma)^{\frac{1}{4}}}{3-\gamma}\tau^{\frac{3-\gamma}{4}}.
\end{equation}
Using \eqref{4.8-1}, \eqref{4.9}-\eqref{4.9-1}, we obtain
\begin{equation}\label{3.21}
	\begin{aligned}
		G^\prime&=(c^{\frac{1}{2}}w_x)^\prime+\frac{\a}{(1+t)^\lambda}h^\prime(\t)-\frac{\lambda\a}{(1+t)^{\lambda+1}}h(\t)\\
		%&=(c^{\frac{1}{2}}w_x)'+\frac{\a}{(1+t)^\lambda}c^{\frac{1}{2}}\t'-\frac{\lambda\a}{(1+t)^{\lambda+1}}h(\t)\\
		&=-\frac{\a c^{\frac{1}{2}}}{2(1+t)^\lambda}({z}+{w})_x-c^{1\over 2}c_w w_x^2
		+\frac{\a}{(1+t)^\lambda}c^{\frac{1}{2}}z_x-\frac{\lambda\a}{(1+t)^{\lambda+1}}h(\t)\\
		&=\frac{\a }{2(1+t)^\lambda}(H-G)-c^{1\over 2}c_w w_x^2-{\lambda\alpha\over {(1+t)^{\lambda+1}}}h(\tau)\\
		&:=-\frac{\a }{2(1+t)^\lambda}( G-H)+R_1,
	\end{aligned}
\end{equation}
and
\begin{equation}\label{3.22}
	\begin{aligned}
		H^\backprime&=(c^{\frac{1}{2}}z_x)^\backprime+\frac{\a}{(1+t)^\lambda}h^\backprime(\t)-\frac{\lambda\a}{(1+t)^{\lambda+1}}h(\t)\\
		&=-\frac{\a c^{\frac{1}{2}}}{2(1+t)^\lambda}({z}+{w})_x-c^{1\over 2}c_w z_x^2
		+\frac{\a}{(1+t)^\lambda}c^{\frac{1}{2}}w_x-\frac{\lambda\a}{(1+t)^{\lambda+1}}h(\t)\\
		&=\frac{\a }{2(1+t)^\lambda}(G-H)-c^{1\over 2}c_w z_x^2-{\lambda\alpha\over {(1+t)^{\lambda+1}}}h(\tau)\\
		&:=-\frac{\a }{2(1+t)^\lambda}(H-G)+R_2.
	\end{aligned}
\end{equation}
Noticing from \eqref{1.4-1} and \eqref{1.5-1} that $c_{w}=\frac{\gamma+1}{2(\gamma-1)}K_{c}\phi^{\frac{2}{\gamma-1}}\geq 0$, then
$$
R_1:=-c^{1\over 2}c_w w_x^2-{\alpha\lambda\over {(1+t)^{\lambda+1}}}h(\tau)\leq0,
$$
$$
R_2:=-c^{1\over 2}c_w z_x^2-{\alpha\lambda\over {(1+t)^{\lambda+1}}}h(\tau)\leq0,
$$
which together with \eqref{3.21}-\eqref{3.22}, yields
$$
G^\prime\leq-\frac{\a }{2(1+t)^\lambda}( G-H),
$$
$$
H^\backprime\leq-\frac{\a }{2(1+t)^\lambda}(H-G).
$$

It follows from \eqref{1.4}, \eqref{1.5-1}, \eqref{1.10} and \eqref{b.1}-\eqref{b.02} that there exists a positive constant $M$ depending only on $C_{0}, K, \gamma, \lambda$ and $\alpha$ such that $\sup_{x\in \R}|(G(0,x), H(0,x))|\leq M$. Now let
\begin{equation}\label{3.20-1}
G_1(t,x)=G(t,x)-M,\quad H_1(t,x)=H(t,x)-M,
\end{equation}
then $G_1(0,x)\leq 0, H_1(0,x)\leq 0$ and
\begin{equation}\label{3.20-2}
\begin{aligned}
&G_1^\prime={G}^\prime\leq-\frac{\a }{2(1+t)^\lambda}( G-H)=-\frac{\a }{2(1+t)^\lambda}G_1
+\frac{\alpha}{2(1+t)^\lambda}H_1,\\
&H_1^\backprime={H}^\backprime\leq-\frac{\a }{2(1+t)^\lambda}(H-G)=-\frac{\a }{2(1+t)^\lambda}H_1
+\frac{\alpha}{2(1+t)^\lambda}G_1.
\end{aligned}
\end{equation}
Since $\lambda \neq 1$, we denote $A(t)=e^{\frac{\alpha(1+t)^{1-\lambda}}{2(1-\lambda)}}$ and
$$
G_2(t,x)=A(t)G_1(t,x),\quad H_2(t,x)=A(t)H_1(t,x),
$$
then $G_2(0,x)\leq 0, H_{2}(0,x)\leq 0$ and
\begin{equation}\label{4.20}
\begin{aligned}
&G_2^\prime=A(t)G_{1}^{\prime}+\frac{\alpha}{2(1+t)^{\lambda}}A(t)G_1\leq \frac{\alpha}{2(1+t)^\lambda}A(t)H_1=\frac{\alpha}{2(1+t)^\lambda}H_2,\\
&H_2^\backprime=A(t)H_{1}^{\backprime}+\frac{\alpha}{2(1+t)^{\lambda}}A(t)H_1\leq \frac{\alpha}{2(1+t)^\lambda}A(t)G_1=\frac{\alpha}{2(1+t)^\lambda}G_2.
\end{aligned}
\end{equation}
Using a continuity argument, it is directly from \eqref{4.20} to see
$$
G_2(x_{+}(t),t)\leq0,\quad
H_2(x_{-}(t),t)\leq 0,\quad \text{for all }t\geq 0.
$$
Therefore,
$$
G_1(x_{+}(t),t)\leq0,\quad H_1(x_{-}(t),t)\leq0,\quad \text{for all }t\geq 0,
$$
which together with \eqref{3.20-1} yields that
\begin{equation}\label{04.21}
	G(t,x)\leq M,\quad H(t,x)\leq M,\quad \text{for all }(t,x)\in \R_{+}\times \R.
\end{equation}
Hence by \eqref{b.1}, \eqref{b.2} and \eqref{04.21}, we obtain
$$
G(t,x)+H(t,x)=c^{1\over 2}(w+z)_x+{2\alpha\over {(1+t)^\lambda}}h(\tau)\leq 2M,
$$
which together with $\eqref{1.1}_{1}$ implies that
\begin{equation}\label{521}
	c^{1\over 2}\tau_t+{\alpha\over {(1+t)^\lambda}}h(\tau)\leq M.
\end{equation}

Combing \eqref{1.4-1}, \eqref{b.02} with \eqref{521}, we have
$$
(K\gamma)^{1\over 4}\tau^{-\frac{\gamma+1}{4}}\tau_t+{\alpha\over {(1+t)^\lambda}}\frac{4(K\gamma)^{\frac{1}{4}}}{3-\gamma}\tau^{\frac{3-\gamma}{4}}\leq M,
$$
%which implies that
%\begin{equation}\label{04.23}
	%\frac{4}{3-\gamma}(\tau^{\frac{3-\gamma}{4}})_t+{\alpha\over %{(1+t)^\lambda}}\frac{4}{3-\gamma}\tau^{\frac{3-\gamma}{4}}\leq %(K\gamma)^{-\frac{1}{4}}M.
%\end{equation}
%By (\ref{04.23}), we have
%$$
%(\tau^{\frac{3-\gamma}{4}})_t+{\alpha\over {(1+t)^\lambda}}\big(\tau^{\frac{3-\gamma}{4}}-\tilde{C}_0^{\frac{\gamma-3}{4}}\big)
%\leq \frac{3-\gamma}{4}(K\gamma)^{-\frac{1}{4}}M,
%$$
%then
%$$
%(\tau^{\frac{3-\gamma}{4}})_t+{\alpha\over {(1+t)^\lambda}}\tau^{\frac{3-\gamma}{4}}
%\leq \frac{3-\gamma}{4}(K\gamma)^{-\frac{1}{4}}M,
%$$
which implies that
\begin{equation}\label{04.25}
	\left(e^{\frac{\alpha(1+t)^{1-\lambda}}{1-\lambda}}\tau^{\frac{3-\gamma}{4}}\right)_t\leq \left(\frac{3-\gamma}{4}(K\gamma)^{-\frac{1}{4}}M\right)e^{\frac{\alpha(1+t)^{1-\lambda}}{1-\lambda}}.
\end{equation}
Integrating (\ref{04.25}) from $0$ to $t$, one obtains
\begin{equation}\label{04.25-1}
\begin{aligned}
	e^{\frac{\alpha(1+t)^{1-\lambda}}{1-\lambda}}\tau^{\frac{3-\gamma}{4}}
	-e^{\frac{\alpha}{1-\lambda}}\tau_0^{\frac{3-\gamma}{4}}
	&\leq \frac{3-\gamma}{4}(K\gamma)^{-\frac{1}{4}}M\int_0^t e^{\frac{\alpha(1+s)^{1-\lambda}}{1-\lambda}}ds\\
	&\leq \frac{3-\gamma}{4}(K\gamma)^{-\frac{1}{4}}Me^{\frac{\alpha(1+t)^{1-\lambda}}{1-\lambda}}t,
\end{aligned}
\end{equation}
therefore
\begin{equation}\label{04.25-2}
\begin{aligned}
	\tau^{\frac{3-\gamma}{4}}
	&\leq \left[e^{\frac{\alpha}{1-\lambda}}\tau_0^{\frac{3-\gamma}{4}}
	+\frac{3-\gamma}{4}(K\gamma)^{-\frac{1}{4}}M e^{\frac{\alpha(1+t)^{1-\lambda}}{1-\lambda}}t\right]e^{-\frac{\alpha(1+t)^{1-\lambda}}{1-\lambda}}\\
	&= e^{\frac{\alpha(1-(1+t)^{1-\lambda})}{1-\lambda}}\tau_0^{\frac{3-\gamma}{4}}+\frac{3-\gamma}{4}(K\gamma)^{-\frac{1}{4}}M t,
\end{aligned}
\end{equation}
i.e.,
\begin{equation}\label{04.25-3}
\begin{aligned}
	\rho &\geq \left[\tau_0^{\frac{3-\gamma}{4}}+\frac{3-\gamma}{4}(K\gamma)^{-\frac{1}{4}}M t\right]^{\frac{4}{\gamma-3}}\\
	&\geq \left[C_0^{\frac{3-\gamma}{4}}+\frac{3-\gamma}{4}(K\gamma)^{-\frac{1}{4}}M t\right]^{\frac{4}{\gamma-3}}\\
	&\geq C(1+t)^{\frac{4}{\gamma-3}},
\end{aligned}
\end{equation}
where we have used $1<\gamma<3$, \eqref{1.10} and
\begin{equation}\label{04.25-4}
C:=\left(\max\{C_0^{\frac{3-\gamma}{4}}, \frac{3-\gamma}{4}(K\gamma)^{-\frac{1}{4}}M\}\right)^{\frac{4}{\gamma-3}}.
\end{equation}
%\begin{equation}\label{e.3}
	%M_1=\tilde{C}_0^{{3-\gamma}\over 4},\ \  %M_2=\frac{3-\gamma}{4}(K\gamma)^{-\frac{1}{4}}M+\alpha %\tilde{C}_0^{\frac{\gamma-3}{4}}.
%\end{equation}
\\

\textbf{Step 3.} For $\lambda=1$, we denote $A_1(t)=(1+t)^{\frac{\alpha}{2}}$ and
$$
G_3(t,x)=A_1(t)G_{1}(t,x),\quad H_{3}(t,x)=A_1(t)H_1(t,x),
$$
then it follows from \eqref{3.20-2} that $G_3(0,x)\leq 0, H_3(0,x)\leq 0$ and
\begin{equation}\label{4.20-1}
\begin{aligned}
&G_3^\prime=A_1(t)G_{1}^{\prime}+\frac{\alpha}{2(1+t)}A_1(t)G_1\leq \frac{\alpha}{2(1+t)}A_1(t)H_1=\frac{\alpha}{2(1+t)}H_3,\\
&H_3^\backprime=A_1(t)H_{1}^{\backprime}+\frac{\alpha}{2(1+t)}A_1(t)H_1\leq \frac{\alpha}{2(1+t)}A_1(t)G_1=\frac{\alpha}{2(1+t)}G_3.
\end{aligned}
\end{equation}
Hence, by a similar argument in Step 2, we obtain
$$
G(t,x)\leq M,\quad H(t,x)\leq M,\quad \text{for all }(t,x)\in \R_{+}\times \R,
$$
which implies that
\begin{equation}\label{4.20-2}
G(t,x)+H(t,x)=c^{\frac{1}{2}}(w+z)_{x}+\frac{2\alpha}{(1+t)}h(\t)\leq 2M.
\end{equation}
Using $\eqref{1.1}_1$, \eqref{1.4-1} and \eqref{b.02}, then \eqref{4.20-2} becomes
\begin{equation}\label{4.20-3}
	(\tau^{\frac{3-\gamma}{4}})_t+{\alpha\over {(1+t)}}\tau^{\frac{3-\gamma}{4}}
	\leq \frac{3-\gamma}{4}(K\gamma)^{-\frac{1}{4}}M,
\end{equation}
and a direct calculation shows that
\begin{equation}\label{4.20-4}
	\big((1+t)^{\a}\tau^{\frac{3-\gamma}{4}}\big)_t
	\leq \frac{3-\gamma}{4}(K\gamma)^{-\frac{1}{4}}M(1+t)^{\a}.
\end{equation}
Integrating \eqref{4.20-4} from $0$ to $t$, we have
$$
\begin{aligned}
(1+t)^{\a}\tau^{\frac{3-\gamma}{4}}-\tau_{0}^{\frac{3-\gamma}{4}}&\leq \frac{3-\gamma}{4}(K\gamma)^{-\frac{1}{4}}M\int_{0}^{t}(1+s)^{\a}ds
\\&\leq \frac{3-\gamma}{4}(K\gamma)^{-\frac{1}{4}}M(1+t)^{\a}t,
\end{aligned}
$$
and therefore
\begin{equation}\label{04.25-5}
\begin{aligned}
	\rho &\geq \left[(1+t)^{-\a}\tau_0^{\frac{3-\gamma}{4}}+\frac{3-\gamma}{4}(K\gamma)^{-\frac{1}{4}}M t\right]^{\frac{4}{\gamma-3}}\\
	&\geq \left[\tau_0^{\frac{3-\gamma}{4}}+\frac{3-\gamma}{4}(K\gamma)^{-\frac{1}{4}}M t\right]^{\frac{4}{\gamma-3}}\\
	&\geq C(1+t)^{\frac{4}{\gamma-3}},
\end{aligned}
\end{equation}
where we have used $1<\gamma<3$, \eqref{1.10} and \eqref{04.25-4}.\\

\textbf{Step 4.} For $\lambda<0$, we introduce two new variables
\begin{equation}\label{4.1}
	\hat{G}(t,x)=c^{1\over 2}w_x+{\alpha-\lambda\over {(1+t)^\lambda}}h(\tau),
\end{equation}
\begin{equation}\label{4.2}
	\hat{H}(t,x)=c^{1\over 2}z_x+{\alpha-\lambda\over {(1+t)^\lambda}}h(\tau),
\end{equation}
where $h(\tau)$ is defined \eqref{b.02}. By a similar calculation in \eqref{3.21}-\eqref{3.22}, one has
\begin{equation}\label{4.11}
	\begin{aligned}
	 \hat{G}'&=\frac{\alpha-2\lambda}{2(1+t)^\lambda} \hat{H}
		-\frac{\a }{2(1+t)^\lambda}\hat{G}-c^{1\over 2}c_w w_x^2
		+\frac{\lambda(\alpha-\lambda)(1-(1+t)^{\lambda-1})}{(1+t)^{2\lambda}}h(\tau),\\
		&:=\frac{\alpha-2\lambda}{2(1+t)^\lambda} \hat{H}
		-\frac{\a }{2(1+t)^\lambda}\hat{G}+V_1,\\
		%&=(c^{\frac{1}{2}}w_x)'+\frac{\alpha-\lambda}{(1+t)^\lambda}h'(\t)-\frac{\lambda(\alpha-\lambda)}{(1+t)^{\lambda+1}}h(\t)\\
		%&=(c^{\frac{1}{2}}w_x)'+\frac{\alpha-\lambda}{(1+t)^\lambda}c^{\frac{1}{2}}\t'-\frac{\lambda(\alpha-\lambda)}{(1+t)^{\lambda+1}}h(\t)\\
		%&=-\frac{\a c^{\frac{1}{2}}}{2(1+t)^\lambda}({z}+{w})_x-c^{1\over 2}c_w w_x^2
		%+\frac{\alpha-\lambda}{(1+t)^\lambda}c^{\frac{1}{2}}z_x-\frac{\lambda(\alpha-\lambda)}{(1+t)^{\lambda+1}}h(\t)\\
		%&=-\frac{\a }{2(1+t)^\lambda}\left(\hat{H}+\hat{G}-2{\alpha-\lambda\over {(1+t)^\lambda}}h(\tau)\right)
		%+\frac{\alpha-\lambda}{(1+t)^\lambda}\left(\hat{H}-{\alpha-\lambda\over {(1+t)^\lambda}}h(\tau)\right)\\&\ \
		%%-c^{1\over 2}c_w w_x^2-\lambda{\alpha-\lambda\over {(1+t)^{\lambda+1}}}h(\tau)\\
		\hat{H}^\backprime
		&=-\frac{\a }{2(1+t)^\lambda}\hat{H}+\frac{\alpha-2\lambda}{2(1+t)^\lambda} \hat{G}
		-c^{1\over 2}c_w z_x^2
		+\frac{\lambda(\alpha-\lambda)(1-(1+t)^{\lambda-1})}{(1+t)^{2\lambda}}h(\tau)\\
		&:=-\frac{\a }{2(1+t)^\lambda}\hat{H}+\frac{\alpha-2\lambda}{2(1+t)^\lambda} \hat{G}+V_2,
	\end{aligned}
\end{equation}
where
$$
V_1:=-c^{1\over 2}c_w w_x^2
+\frac{\lambda(\alpha-\lambda)(1-(1+t)^{\lambda-1})}{(1+t)^{2\lambda}}h(\tau)\leq0,
$$
$$
V_2:=-c^{1\over 2}c_w z_x^2
+\frac{\lambda(\alpha-\lambda)(1-(1+t)^{\lambda-1})}{(1+t)^{2\lambda}}h(\tau)\leq0,$$
for $\lambda <0$.

%Similarly,
%\begin{equation}\label{4.12}
	%\begin{aligned}
		%\hat{H}^\backprime
		%=-\frac{\a %}{2(1+t)^\lambda}\hat{H}+\frac{\alpha-2\lambda}{2(1+t)^\lambda} \hat{G}
		%-c^{1\over 2}c_w z_x^2
		%+\frac{\lambda(\alpha-\lambda)(1-(1+t)^{\lambda-1})}{(1+t)^{2\lambda}}h(\%tau).
	%\end{aligned}
%\end{equation}
%Then (\ref{4.11}) and (\ref{4.12}) turn into
%\begin{equation}\label{4.13}
%	\begin{aligned}
%		\hat{G}^\prime= \frac{\alpha-2\lambda}{2(1+t)^\lambda} \hat{H}
%		-\frac{\a }{2(1+t)^\lambda}\hat{G}+V_1,
%	\end{aligned}
%\end{equation}
%\begin{equation}\label{4.14}
%	\begin{aligned}
%		\hat{H}^\backprime= \frac{\alpha-2\lambda}{2(1+t)^\lambda} \hat{G}
%		-\frac{\a }{2(1+t)^\lambda}\hat{H}+V_2,
%	\end{aligned}
%\end{equation}
%where $$V_1=-c^{1\over 2}c_w w_x^2
%-\frac{\lambda(\alpha-\lambda)(1-(1+t)^{\lambda-1})}{(1+t)^{2\lambda}}h(\tau)\l%eq0,$$
%$$V_2=-c^{1\over 2}c_w z_x^2
%-\frac{\lambda(\alpha-\lambda)(1-(1+t)^{\lambda-1})}{(1+t)^{2\lambda}}h(\tau)\l%eq0,$$
%for $\lambda \leq0$.
Let
$$
\hat{G}_1(t,x)=\hat{G}(t,x)-e^{-\frac{\lambda(1+t)^{1-\lambda}}{1-\lambda}}M,\quad
\hat{H}_1(t,x)=\hat{H}(t,x)-e^{-\frac{\lambda(1+t)^{1-\lambda}}{1-\lambda}}M,
$$
where $M$, like the constant $M$ defined in \eqref{3.20-1}, is a big enough constant determined by the initial data, such that
\begin{equation}\label{4.11-1}
\hat{G}_1(0,x)=\hat{G}(0)-e^{-\frac{\lambda}{1-\lambda}}M\leq 0,\quad
\hat{H}_1(0,x)=\hat{H}(0)-e^{-\frac{\lambda}{1-\lambda}}M\leq 0.
\end{equation}
A direct calculation by using \eqref{4.11} shows that
\begin{equation}\label{4.15}
	\begin{aligned}
		\hat{G}_1^\prime&=\hat{G}^\prime+\lambda(1+t)^{-\lambda}e^{-\frac{\lambda(1+t)^{1-\lambda}}{1-\lambda}}M\\
		&\leq \frac{\alpha-2\lambda}{2(1+t)^\lambda}\left(\hat{H}_1+e^{-\frac{\lambda(1+t)^{1-\lambda}}{1-\lambda}}M\right)-\frac{\a }{2(1+t)^\lambda}\left(\hat{G}_1+e^{-\frac{\lambda(1+t)^{1-\lambda}}{1-\lambda}}M\right)+\frac{\lambda e^{-\frac{\lambda(1+t)^{1-\lambda}}{1-\lambda}}M}{(1+t)^{\lambda}}\\
		&= -\frac{\a }{2(1+t)^\lambda}\hat{G}_1
		+\frac{\alpha-2\lambda}{2(1+t)^\lambda}\hat{H}_1
		+\left(\frac{\alpha-2\lambda}{2(1+t)^\lambda}-\frac{\a }{2(1+t)^\lambda}+\frac{\lambda}{(1+t)^\lambda}\right)e^{-\frac{\lambda(1+t)^{1-\lambda}}{1-\lambda}}M\\
		&=-\frac{\a }{2(1+t)^\lambda}\hat{G}_1
		+\frac{\alpha-2\lambda}{2(1+t)^\lambda}\hat{H}_1,
	\end{aligned}
\end{equation}
and similarly
\begin{equation}\label{4.16}
	\begin{aligned}
		\hat{H}_1^\backprime
		&\leq-\frac{\a }{2(1+t)^\lambda}\hat{H}_1
		+\frac{\alpha-2\lambda}{2(1+t)^\lambda}\hat{G}_1.
	\end{aligned}
\end{equation}
Since $\lambda< 0$, denoting $A(t)=e^{\frac{\alpha(1+t)^{1-\lambda}}{2(1-\lambda)}}$ and
$$
\hat{G}_2(t,x)=A(t)\hat{G}_1(t,x),\quad
\hat{H}_2(t,x)=A(t)\hat{H}_1(t,x),
$$
then it follows from \eqref{4.11}-\eqref{4.11-1}, \eqref{4.15}-\eqref{4.16} that
$$
\hat{G}_2(0,x)\leq0,\quad \hat{H}_2(0,x)\leq0,
$$
and
$$
\hat{G}_2^\prime\leq \frac{\alpha-2\lambda}{2(1+t)^\lambda}\hat{H}_2,\quad
\hat{H}_2^\backprime\leq \frac{\alpha-2\lambda}{2(1+t)^\lambda}\hat{G}_2.
$$
Hence by a similar argument in Step 2, we obtain
%$$
%\hat{G}_2(x(t),t)\leq0,
%\ \
%\hat{H}_2(x(t),t)\leq0,
%$$
%therefore,
%$$
%\hat{G}_1(x(t),t)\leq0,
%\ \
%\hat{H}_1(x(t),t)\leq0,
%$$
%since $$\hat{G}_2=A(t)\hat{G}_1\leq A(0)\hat{G}_1(0)=A(0)(\hat{G}(0)-e^{-\frac{\lambda}{1-\lambda}}M)\leq0,$$
%$$\hat{H}_2=A(t)\hat{H}_1\leq A(0)\hat{H}_1(0)=A(0)(\hat{H}(0)-e^{-\frac{\lambda}{1-\lambda}}M)\leq0,$$
%thus
%\begin{equation}\label{4.21}
%	\hat{G}\leq e^{-\frac{\lambda(1+t)^{1-\lambda}}{1-\lambda}}M,
%\end{equation}
%\begin{equation}\label{4.22}
%	\hat{H}\leq e^{-\frac{\lambda(1+t)^{1-\lambda}}{1-\lambda}}M.
%\end{equation}
%And by (\ref{4.1}) (\ref{4.2}) (\ref{4.21}) and (\ref{4.22})
$$
\hat{G}(t,x)+\hat{H}(t,x)=c^{1\over 2}(w+z)_x+{2(\alpha-\lambda)\over {(1+t)^\lambda}}h(\tau)\leq 2e^{-\frac{\lambda(1+t)^{1-\lambda}}{1-\lambda}}M,
$$
which together with $\eqref{1.1}_1$, yields
\begin{equation}\label{520}
	c^{1\over 2}\tau_t+{\alpha-\lambda\over {(1+t)^\lambda}}h(\tau)\leq e^{-\frac{\lambda(1+t)^{1-\lambda}}{1-\lambda}}M.
\end{equation}

Thus using a similar calculation in \eqref{521}-\eqref{04.25-3}, we have
$$
\begin{aligned}
	\tau^{\frac{3-\gamma}{4}}
	&\leq \left[e^{\frac{\alpha-\lambda}{1-\lambda}}\tau_0^{\frac{3-\gamma}{4}}+\frac{3-\gamma}{4}(K\gamma)^{-\frac{1}{4}}M e^{\frac{(\alpha-2\lambda)(1+t)^{1-\lambda}}{1-\lambda}}t\right]
	e^{-\frac{(\alpha-\lambda)(1+t)^{1-\lambda}}{1-\lambda}}\\
	&=e^{\frac{(\alpha-\lambda)(1-(1+t)^{1-\lambda})}{1-\lambda}}\tau_0^{\frac{3-\gamma}{4}}+\frac{3-\gamma}{4}(K\gamma)^{-\frac{1}{4}}M e^{\frac{-\lambda(1+t)^{1-\lambda}}{1-\lambda}}t,
\end{aligned}
$$
i.e.,
\begin{equation}\label{04.25-6}
\begin{aligned}
\rho&\geq \left[e^{\frac{(\alpha-\lambda)(1-(1+t)^{1-\lambda})}{1-\lambda}}\tau_0^{\frac{3-\gamma}{4}}+\frac{3-\gamma}{4}(K\gamma)^{-\frac{1}{4}}M e^{\frac{-\lambda(1+t)^{1-\lambda}}{1-\lambda}}t\right]^{\frac{4}{\gamma-3}}
\\&\geq \left[\tau_0^{\frac{3-\gamma}{4}}+\frac{3-\gamma}{4}(K\gamma)^{-\frac{1}{4}}M e^{\frac{-\lambda(1+t)^{1-\lambda}}{1-\lambda}}t\right]^{\frac{4}{\gamma-3}}\\
&\geq C\left[1+e^{\frac{-\lambda(1+t)^{1-\lambda}}{1-\lambda}}t\right]^{\frac{4}{\gamma-3}},
\end{aligned}
\end{equation}
where we have used $1<\gamma<3$, \eqref{1.10} and \eqref{04.25-4}.
\\

\textbf{Step 5.} Now we are going to prove the blow-up of $\|(\tau_{x},u_{x})\|_{L^{\infty}}$. Since the case $\lambda\geq \frac{\alpha(\gamma-1)}{\gamma-3}$ is solved in Theorem \ref{thm1.2}, then we only need to prove the rest case $\lambda< \frac{\alpha(\gamma-1)}{\gamma-3}$.
%%With the help of lower bound estimates on the density obtained in Step 1-Step 3, the main proof idea is similar with the case
%%$$
%%\gamma>3, \ \ 0<\frac{\alpha(\gamma-1)}{\gamma-3}< \lambda <1
%%$$
%%in Theorem \ref{thm1.3} \\
Noticing that $\lambda<\frac{\alpha(\gamma-1)}{\gamma-3}< 0$ for $1<\gamma<3$, it follows from \eqref{1.15} that there exists $t_0>0$ such that
$$
(1+t_{0})^{\lambda-1}=\frac{\alpha(\gamma-1)}{\lambda(\gamma-3)}
$$
and
$$
a_0\geq 0,\ \ \forall\ 0<t\leq t_0;\ \ a_0<0,\ \ \forall\ t>t_0.
$$
Using \eqref{1.4}, \eqref{1.16}-\eqref{1.15} and Lemma \ref{lem1.1}, there exist two constants $\tilde{K}_3$ and $\tilde{K}_4$, depending only on $C_{0}, \gamma, K, \lambda$ and $\alpha$, such that
\begin{equation}\label{4.26}
0\leq a_0\leq \tilde{K}_3, \quad a_2\geq \tilde{K}_4>0,\quad \forall\ 0<t\leq t_0.
\end{equation}
When $t\geq t_0$, using \eqref{04.25-6} and a direct calculation shows that
\begin{equation}\label{4.26-1}
	\begin{aligned}
		&\int_{t_0}^{\infty}a_2(x,s)ds\\
		= &\frac{K_{c}(\gamma+1)}{2(\gamma-1)}\left(\frac{2\sqrt{K\gamma}}{\gamma-1}\right)^{-\frac{\gamma-3}{2(\gamma-1)}}
		\int_{t_0}^{+\infty}\rho^{\frac{3-\gamma}{4}}e^{-\frac{\alpha(3\gamma-1)}{2(\gamma-3)(1-\lambda)}(1+s)^{1-\lambda}}ds\\
		\geq & \frac{K_{c}(\gamma+1)}{2(\gamma-1)}\left(\frac{2\sqrt{K\gamma}}{\gamma-1}\right)^{-\frac{\gamma-3}{2(\gamma-1)}}C^{\frac{3-\gamma}{4}}
		\int_{t_0}^{+\infty} \left[1+ e^{\frac{-\lambda(1+s)^{1-\lambda}}{1-\lambda}}s\right]^{-1}e^{-\frac{\alpha(3\gamma-1)}{2(\gamma-3)(1-\lambda)}(1+s)^{1-\lambda}}ds\\
		\geq&
		{\left\{\begin{matrix}
				\frac{K_{c}(\gamma+1)}{2(\gamma-1)}\left(\frac{2\sqrt{K\gamma}}{\gamma-1}\right)^{-\frac{\gamma-3}{2(\gamma-1)}}C^{\frac{3-\gamma}{4}}
                \disp\int_{t_0}^{+\infty}\left[e^{\frac{\lambda(1+s)^{1-\lambda}}{1-\lambda}}+ s\right]^{-1}e^{\left(\frac{\alpha( 3\gamma-1)}{2(3-\gamma)}+\lambda\right)\frac{(1+s)^{1-\lambda}}{1-\lambda}}ds \\
                \text{for } \frac{\alpha( 3\gamma-1)}{2(3-\gamma)}+\lambda <0
                \quad\quad\quad\quad\quad\quad\quad\quad\quad\quad\quad\quad\quad\quad\quad\quad\quad\quad\quad\quad\quad\quad\quad\quad
				\\\frac{K_{c}(\gamma+1)}{2(\gamma-1)}\left(\frac{2\sqrt{K\gamma}}{\gamma-1}\right)^{-\frac{\gamma-3}{2(\gamma-1)}}C^{\frac{3-\gamma}{4}}
                e^{\left(\frac{\alpha( 3\gamma-1)}{2(3-\gamma)(1-\lambda)}+\frac{\lambda}{1-\lambda}\right)}\disp\int_{t_0}^{+\infty}\left[e^{\frac{\lambda(1+s)^{1-\lambda}}{1-\lambda}}+ s\right]^{-1}ds\\
                \text{for }\frac{\alpha( 3\gamma-1)}{2(3-\gamma)}+\lambda\geq0  \quad\quad\quad\quad\quad\quad\quad\quad\quad\quad\quad\quad\quad\quad\quad\quad\quad\quad\quad\quad\quad\quad\quad\quad
			\end{matrix}
			\right.}\\
		\geq &{\left\{\begin{matrix}
		\tilde{N}_1>0, &\quad \text{for } \frac{\alpha( 3\gamma-1)}{2(3-\gamma)}+\lambda <0\\
		+\infty,&\quad \text{for }\frac{\alpha( 3\gamma-1)}{2(3-\gamma)}+\lambda\geq0.
	\end{matrix}
     \right.}
	\end{aligned}
\end{equation}
%and
%\begin{equation}
	%\begin{aligned}
		%&\int_{t_0}^{+\infty}s^{-1}e^{\frac{\lambda}{1-\lambda}(1+s)^{1-\lambda}}e^{\frac{\alpha(3\gamma-1)}{2(3-\gamma)(1-\lambda)}(1+s)^{1-\lambda}}ds
		%\\ = &\int_{t_0}^{+\infty}s^{-1}
		%e^{\left(\frac{\alpha( 3\gamma-1)}{2(3-\gamma)}+\lambda\right)\frac{(1+s)^{1-\lambda}}{1-\lambda}}ds
		%\\ \geq& {\left\{\begin{matrix}
				%\disp\int_{t_0}^{+\infty}s^{-1}e^{\left(\frac{\alpha( 3\gamma-1)}{2(3-\gamma)}+\lambda\right)\frac{(1+s)^{1-\lambda}}{1-\lambda}}ds:=\tilde{N}_2>0,
				%&\quad \text{for } \frac{\alpha( 3\gamma-1)}{2(3-\gamma)}+\lambda <0,
				%\\ e^{\left(\frac{\alpha( 3\gamma-1)}{2(3-\gamma)(1-\lambda)}+\frac{\lambda}{1-\lambda}\right)}\disp\int_{t_0}^{+\infty}s^{-1}ds =+\infty,
				%&\quad  \text{for }\frac{\alpha( 3\gamma-1)}{2(3-\gamma)}+\lambda\geq0.
			%\end{matrix}
		%	\right.}
	%\end{aligned}
%\end{equation}
Assume there exists one point $x_0$ such that
$$
y(x_0,0)<-\max\big\{\frac{1}{\tilde{N}_1},\sqrt{\frac{\tilde{K}_3}{\tilde{K}_4}}\big\}:=-N_1,
$$
which means (\ref{1.12}) holds with
$$
K_3:=\frac{\alpha(\gamma-1)}{K_{c}(\gamma-3)},\quad  K_4:=N_1e^{-\frac{\alpha(3\gamma-1)}{2(\gamma-3)(1-\lambda)}}.
$$
Noticing from \eqref{1.14} and \eqref{4.26}, along the forward characteristic starting form $x_0$, we have
\begin{equation}\label{30}
y^\prime<0\quad \text{ and }\quad y(x_{+}(t),t)<-N_1,\quad \forall\ 0<t\leq t_0.
\end{equation}
If $y(x_{+}(t),t)$ does not blow up before $t_0$, then
$$
y^\prime=a_0-a_{2}y^2<-a_2y^2,\quad \forall\ t\geq t_0,
$$
which implies that
\begin{equation}\label{31}
0>y^{-1}(x_{+}(t),t)\geq y^{-1}(x_{+}(t_0),t_0)+\int_{t_0}^{t}a_2(x_{+}(s),s)ds.
\end{equation}
Since it follows from \eqref{30} that $-\frac{1}{N_1}<y^{-1}(x_{+}(t_0),t_0)<0$,  then using \eqref{4.26-1} and \eqref{31}, we obtain that $y(x_{+}(t),t)$ must blow up in finite time. This holds for $q$ if \eqref{1.13} holds.

On the other hand, it follows from \eqref{1.4} and \eqref{04.25-6} that $\phi^{\frac{\gamma-3}{2(\gamma-1)}}$ will remain bounded for any finite time, then recalling \eqref{1.14-1}-\eqref{1.14-2}, we obtain that $y(x_{+}(t),t)$ or $q(x_{-}(t),t)$ blows up in finite time means $\|(\tau_{x},u_{x})\|_{L^{\infty}}$ blows up in finite time. Therefore the proof of Theorem \ref{thm1.4} is completed.$\hfill\square$

\subsection{Proof of Theorem \ref{thm1.3}}
 \ \

We divide the proof into two cases: $\lambda>1$ and $\lambda<1$.\\

\textbf{Case 1.} For $\gamma>3$ and $\lambda>1$, it follows from Lemma \ref{lem1.1} and \eqref{1.16}-\eqref{1.15} that there exist two positive constants $\tilde{K}_5,\tilde{K}_6$, which depends only on $C_0, K, \gamma, \lambda$ and $\alpha$, such that
$$
|a_0|\leq \tilde{K}_5,\quad a_2\geq \tilde{K}_6,\quad \text{for all }t\geq 0.
$$
Assume there exists one point $x_0$ such that
\begin{equation}\label{2.1}
	y(x_0,0)<-\sqrt{\frac{\tilde{K}_5}{\tilde{K}_6}}:=-N_2,
\end{equation}
which means \eqref{1.8} holds with

$$
K_5:=\frac{\alpha(\gamma-1)}{K_{c}(\gamma-3)},\quad  K_6:=N_2e^{-\frac{\alpha(3\gamma-1)}{2(\gamma-3)(1-\lambda)}}.
$$
Noticing from \eqref{2.1}, there exists a $\var>0$ such that
$$
y(x_0,0)<-(1+\var)N_2.
$$
Recalling \eqref{1.14}, along the forward characteristic lines $x_{+}(t)$ starting from $x_0$, we have
$$
  y^\prime(x_{+}(t),t)<0\quad\text{ and }\quad y(x_{+}(t),t)<-(1+\var)N_2,\quad \forall \ t\geq 0,
$$
which together with the definition of $N_2$ implies that
$$
a_0-a_2\frac{y^2(x_{+}(t),t)}{(1+\var)^2}<0,\quad \forall\ t\geq 0.
$$
Hence, using \eqref{1.14} and $a_2>0$, we have
$$
  y^\prime(x_{+}(t),t)<\left(-1+\frac{1}{(1+\var)^2}\right)a_2y^2,\quad \forall\ t\geq 0.
$$
Integrating both sides of the above inequality with respect to $t$ along the forward characteristic $x_{+}(t)$ starting form $x_0$, we obtain
\begin{equation}\label{2.2}
0>y^{-1}(x_{+}(t),t)\geq y^{-1}(x_0,0)+\left(1-\frac{1}{(1+\var)^2}\right)\int_{0}^{t}a_2(x_{+}(s),s)ds.
\end{equation}
Noticing that
$$
\int_{0}^{+\infty}a_2(x_{+}(s),s)ds\geq \int_{0}^{+\infty}\tilde{K}_6ds=+\infty,
$$
which shows that the right hand side of (\ref{2.2}) will be positive in finite time, then $y(x_{+}(t),t)$ must blow up in finite time. Similarly, if \eqref{1.9} holds, we will obtain that $q(x_{-}(t),t)$ blows up in finite time. Moreover, it follows from \eqref{1.4} and Lemma \ref{lem1.1} that $\phi^{\frac{\gamma+1}{2(\gamma-1)}}$ and $\phi^{\frac{\gamma-3}{2(\gamma-1)}}$ remain bounded for all time when $\gamma>3$. Hence $y(x_{+}(t),t)$ or $q(x_{-}(t),t)$ blows up means $\|(\tau_{x},u_{x})\|_{L^{\infty}}$ blows up.\\

\textbf{Case 2.} For $\gamma>3$ and $\lambda<1$, $a_0$ may not be uniformly bounded for all $t>0$ and $a_2$ may not have a positive lower bound. However, if $\lambda\leq \min\{\frac{\alpha(\gamma-1)}{\gamma-3},1\}$ or $\alpha=0$, we will have $a_0\leq 0$ for all $t>0$ and hence using \eqref{1.14}, we have
$$
  y^\prime=a_0-a_2y^2<-a_2y^2.
$$
This is the case in Theorem \ref{thm1.1} and is solved in \cite[Theorem 3.1]{Sui-Yu}.

If $\alpha>0$ and $0<\frac{\alpha(\gamma-1)}{\gamma-3}< \lambda <1$, there exists $t_0>0$ satisfying
$$
(1+t_{0})^{\lambda-1}=\frac{\alpha(\gamma-1)}{\lambda(\gamma-3)}
$$
such that
$$
a_0\geq 0,\ \ \forall\ 0<t\leq t_0\ \ \text{and }\  a_0<0,\ \ \forall\ t>t_0.
$$
When $0<t\leq t_0$, using Lemma \ref{lem1.1}, there exists two constants $\hat{K}_5$ and $\hat{K}_6$, which are dependent only on $C_{0}, \gamma, K, \lambda$ and $\alpha$, such that
$$
0\leq a_0\leq \hat{K}_5, \quad a_2\geq \hat{K}_6,\quad \forall\ 0<t\leq t_0.
$$
When $t\geq t_0$, using Lemma \ref{lem1.1}, a direct calculation shows that
\begin{equation}\label{2.3}
\begin{aligned}
&\int_{t_0}^{+\infty}a_2(x(s),s)ds
\\&=\int_{t_0}^{+\infty}\frac{K_{c}(\gamma+1)}{2(\gamma-1)}\phi^{-\frac{\gamma-3}{2(\gamma-1)}}e^{-\frac{\alpha(3\gamma-1)}{2(\gamma-3)(1-\lambda)}(1+s)^{1-\lambda}}ds
\\&=\int_{t_0}^{+\infty}\frac{K_{c}(\gamma+1)}{2(\gamma-1)}\left(\frac{2\sqrt{K\gamma}}{\gamma-1}\tau^{-\frac{\gamma-1}{2}}\right)^{-\frac{\gamma-3}{2(\gamma-1)}}e^{-\frac{\alpha(3\gamma-1)}{2(\gamma-3)(1-\lambda)}(1+s)^{1-\lambda}}ds
\\&\geq \frac{K_{c}(\gamma+1)}{2(\gamma-1)}\left(\frac{2\sqrt{K\gamma}}{\gamma-1}\right)^{-\frac{\gamma-3}{2(\gamma-1)}}\tilde{C}_{0}^{-\frac{\gamma-3}{4}}\int_{t_0}^{+\infty}e^{-\frac{\alpha(3\gamma-1)}{2(\gamma-3)(1-\lambda)}(1+s)^{1-\lambda}}ds
\\&\geq \frac{K_{c}(\gamma+1)}{2(\gamma-1)}\left(\frac{2\sqrt{K\gamma}}{\gamma-1}\right)^{-\frac{\gamma-3}{2(\gamma-1)}}\tilde{C}_{0}^{-\frac{\gamma-3}{4}}\int_{t_0}^{+\infty}(1+s)^{-\lambda}e^{-\frac{\alpha(3\gamma-1)}{2(\gamma-3)(1-\lambda)}(1+s)^{1-\lambda}}ds
\\&\geq \frac{K_{c}(\gamma+1)}{2(\gamma-1)}\left(\frac{2\sqrt{K\gamma}}{\gamma-1}\right)^{-\frac{\gamma-3}{2(\gamma-1)}}\tilde{C}_{0}^{-\frac{\gamma-3}{4}}\frac{2(\gamma-3)}{\alpha(3\gamma-1)}e^{-\frac{\alpha(3\gamma-1)}{2(\gamma-3)(1-\lambda)}}:=\tilde{N}_3.
\end{aligned}
\end{equation}
If there exists a point $x_0$ such that
$$
y(x_0,0)<-\max\{\frac{1}{\tilde{N}_3},\sqrt{\frac{\hat{K}_{5}}{\hat{K}_{6}}}\}:=-N_3,
$$
which is equivalent to \eqref{1.8} with
\begin{equation}
K_5:=\frac{\alpha(\gamma-1)}{K_{c}(\gamma-3)},\quad  K_6:=N_3e^{-\frac{\alpha(3\gamma-1)}{2(\gamma-3)(1-\lambda)}}.
\end{equation}
Then following the same discussion in \eqref{30}-\eqref{31}, we obtain $y(x_{+}(t),t)$ must blow up in finite time. This holds for $q$ if \eqref{1.9} holds. Hence $\|(\tau_{x},u_{x})\|_{L^{\infty}}$ must blow up in finite time. This completes the proof of Theorem \ref{thm1.3}.

$\hfill\square$

%then along the forward characteristic starting form $x_{0}$, we will have
%\begin{equation}\label{e.1}
 % y^\prime\leq 0 \text{ and } y(x_{+}(t),t)<-N_1\quad \forall 0<t\leq t_0.
%\end{equation}
%Hence if $y(x_{+}(t),t)$ does not blow up before $t_0$, we will have
%\begin{equation}
 % y^\prime=a_0-a_{2}y^2<-a_2y^2,\quad \forall t\geq t_0,
%\end{equation}
%and hence for $t>t_0$,
%\begin{equation}\label{e.2}
%0>y^{-1}(x_{+}(t),t)\geq y^{-1}(x_{+}(t_0),t_0)+\int_{t_0}^{t}a_2(x_{+}(s),s)ds.
%\end{equation}
%Since $-\tilde{N}_1\leq -\frac{1}{N_{1}}<y^{-1}(x_{+}(t_0),t_0)<0$, then by (\ref{2.3}), we know $y(x_{+}(t),t)$ will blows up in finite time. This is also true for $q(x_-(t),t)$ if (\ref{1.9}) holds. $\hfill\square$\\

\begin{remark}\label{rem4.1}
	When $\lambda=1$, the authors in \cite{Sui-Yu} have also proved similar blow up results for $\gamma>3,$ $ \alpha\geq \frac{\gamma-3}{\gamma-1}$ and $1<\gamma<3, \alpha\geq 0$. In fact, when $\gamma>3$ and $\lambda=1$, following the same idea in the proofs of Theorem \ref{thm1.4} and Theorem \ref{thm1.3}, the similar blow up phenomenon can be proved for $0\leq \alpha\leq \frac{2(\gamma-3)}{3\gamma-1}$. However, for the rest case $\alpha\in \left(\frac{2(\gamma-3)}{3\gamma-1}, \frac{\gamma-3}{\gamma-1}\right)$, there is no blow up results yet now. The main reason is that we can not get the upper bound of $|{{b_0}\over {b_1}}|$ in (4.3) and (4.4) of \cite{Sui-Yu}.
\end{remark}

\subsection{Proof of Theorem \ref{thm9.1}}
 \ \
We divide the proof into four steps. In steps 1-3, we first give the proof of the lower bound on density, that is, the proof of \eqref{1.200}. In step 4, we prove the blow-up of $\|(\tau_{x},u_{x})\|_{L^{\infty}}$.

\textbf{Step 1.} For $\gamma=3$, the previous notations given in \eqref{1.4}-\eqref{1.4-1} becomes
%\begin{equation*}\label{3}
	%\begin{aligned}
		%p = K\tau^{-3} \ \ \text{for constants $K>0$}.
	%\end{aligned}
%\end{equation*}
%And
\begin{equation}\label{61}
	\begin{aligned}
		\phi :=\int^\infty_\tau c d\tau=\sqrt{3K}\tau^{-1}>0,
	\end{aligned}
\end{equation}
and
%where the nonlinear Lagrangian sound speed $c$ is
\begin{equation}\label{61-1}
	\begin{aligned}
		c :=\sqrt{-p_\tau} = \sqrt{3K}\tau^{-2}:=K_{c}\phi^2.
	\end{aligned}
\end{equation}
%It follows \eqref{61} that
%\begin{equation*}\label{14}
	%c = \sqrt{-p_\tau} = K_c\phi^{2},
%\end{equation*}
where $K_{c}=1/\sqrt{3K}$.
%\begin{equation*}\label{15}
	%\begin{aligned}
		%K_c :=\frac{1}{\sqrt{3K}}.
	%\end{aligned}
%\end{equation*}
%The Riemann invariants
%\begin{equation}\label{62}
	%\begin{aligned}
		%{w}:=u+\phi \ \ \text{and}\ \  {z}:=u-\phi
	%\end{aligned}
%\end{equation}
%satisfy
%\begin{equation}\label{63}
	%\begin{aligned}
		%{w}^\prime:=-\frac{\a}{2(1+t)^\lambda}({z}+{w})
	%\end{aligned}
%\end{equation}
%and
%\begin{equation}\label{64}
	%\begin{aligned}
		%{z}^\backprime:=-\frac{\a}{2(1+t)^\lambda}({z}+{w})
	%\end{aligned}
%\end{equation}
%respectively.\\

We introduce two gradient variables
\begin{equation}\label{68}
	\begin{aligned}
		y_1:=\phi {A}-\frac{\a}{2K_c(1+t)^\lambda}\ln\phi,
	\end{aligned}
\end{equation}
and
\begin{equation}\label{69}
	\begin{aligned}
		q_1:=\phi {B}-\frac{\a}{2K_c(1+t)^\lambda}\ln\phi,
	\end{aligned}
\end{equation}
where $A$ and $B$ are given in Lemma \ref{L6.1} and satisfy
\begin{equation}\label{65}
\begin{aligned}
		&{A}^\prime= -\frac{\a}{2(1+t)^\lambda}({A}+{B})+K_d({A}{B}-{A}^2),\\
		&{B}^\backprime=-\frac{\a}{2(1+t)^\lambda}({A}+{B})+K_d({A}{B}-{B}^2).
\end{aligned}
\end{equation}
For $\gamma=3$, $K_{d}=K_{c}\phi$ from Lemma \ref{L6.1}.\\
\begin{lemma}\label{L6.2}
	For $C^1$ solutions of (\ref{1.1})-(\ref{1.2}),  $y_1$ and $q_1$ satisfy the following Riccati equations:
	\begin{equation}\label{610}
		\begin{aligned}
			y_1^\prime =-a_2y_1^2-a_1y_1+a_0,
		\end{aligned}
	\end{equation}
	\begin{equation}\label{611}
		\begin{aligned}
			q_1^\backprime =-a_2 q_1^2-a_1q_1+a_0,
		\end{aligned}
	\end{equation}
	where
	\begin{equation}\label{612}
		\begin{aligned}
			&a_2=K_c>0,\\
            &a_1=\frac{\a(1+2\ln\phi)}{2(1+t)^\lambda},\\
			&a_0=\a\ln\phi\frac{2\lambda(1+t)^{\lambda-1}-\a(1+\ln\phi)}{4K_c(1+t)^{2\lambda}}.
		\end{aligned}
	\end{equation}
\end{lemma}

\ \\

\textbf{Proof of Lemma \ref{L6.2}}. It follows from \eqref{61}, \eqref{1.1} and \eqref{1.5-1} that
\begin{equation*}
	\begin{aligned}
		\phi^\prime&=\phi_t+c\phi_x
		=\left(\sqrt{3K}\tau^{-1}\right)_t+c\phi_x\\
		&=-\sqrt{3K}\tau^{-2}\tau_t+c\phi_x=-cu_x+c\phi_x\\
		&=-c({z}+\phi)_x+c\phi_x=-c{z}_x-c\phi_x+c\phi_x\\
		&=-c{B},
	\end{aligned}
\end{equation*}
which implies that
\begin{equation}\label{613}
	\begin{aligned}
		{B}=-\frac{1}{c}\phi^\prime.
	\end{aligned}
\end{equation}
Substituting (\ref{613}) into (\ref{65}),  we get
\begin{equation}\label{613-1}
	\begin{aligned}
		{A}^\prime= -\frac{\a}{2(1+t)^\lambda}({A}-\frac{1}{c}\phi^\prime)+K_d(-\frac{1}{c}\phi^\prime{A}-{A}^2).
	\end{aligned}
\end{equation}
Multiplying \eqref{613-1} by $\phi$ to see
\begin{equation}\label{614}
	\begin{aligned}
		{A}^\prime\phi
		-\frac{\a}{2c(1+t)^\lambda}\phi^\prime\phi
		+\frac{K_d}{c}\phi^\prime{A}\phi
		= -\frac{\a}{2(1+t)^\lambda}{A}\phi
		-K_d{A}^2\phi .
	\end{aligned}
\end{equation}
For the left hand side of \eqref{614}, a direct calculation shows that
\begin{equation}\label{614-1}
\begin{aligned}
	&{A}^\prime\phi
	-\frac{\a}{2c(1+t)^\lambda}\phi^\prime\phi
	+\frac{K_d}{c}\phi^\prime{A}\phi \\
	=&{A}^\prime\phi
	-\frac{\a}{2(1+t)^\lambda K_c\phi^{2}}\phi^\prime\phi
	+\frac{K_c \phi }
	{K_c\phi^{2}}\phi^\prime{A}\phi  \\
	=&{A}^\prime\phi
	-\frac{\a}{2K_c(1+t)^\lambda}\phi^\prime\phi^{-1}
	+  \phi^\prime{A}
	\\
	=&\left({A}\phi-\frac{\a}{2K_c(1+t)^\lambda}\ln\phi\right)^\prime
	-\frac{\lambda \a\ln\phi}{2K_c}(1+t)^{-\lambda-1}.
\end{aligned}
\end{equation}
Recalling \eqref{68}, we have
\begin{equation}\label{615}
	\begin{aligned}
		{A}=y\phi^{-1}+\frac{\a\ln\phi}{2K_c(1+t)^\lambda}\phi^{-1}.
	\end{aligned}
\end{equation}
Then the right hand side of \eqref{614} becomes
\begin{equation}\label{615-1}
	\begin{aligned}
		&-\frac{\a}{2(1+t)^\lambda}{A}\phi
		-K_d{A}^2\phi \\
		=&-\frac{\a}{2(1+t)^\lambda}\phi
		\left(y_1\phi^{-1}+\frac{\a\ln\phi}{2K_c(1+t)^\lambda}\phi^{-1}\right)
		-K_c \phi^2\left(y_1\phi^{-1}+\frac{\a\ln\phi}{2K_c(1+t)^\lambda}\phi^{-1}\right)^2\\
		=&-\frac{\a}{2(1+t)^\lambda}y_1
		-\frac{\a^2\ln\phi}{4K_c(1+t)^{2\lambda}}- K_cy_1^2
		-\frac{\a^2(\ln\phi)^2}{4K_c(1+t)^{2\lambda}}
		-\frac{\a\ln\phi}{ (1+t)^\lambda}y\\
		=&-K_cy_1^2
		-\frac{\a+2\a\ln\phi}{2(1+t)^\lambda}y_1
		-\frac{\a^2\ln\phi+\a^2(\ln\phi)^2}{4K_c(1+t)^{2\lambda}}.
	\end{aligned}
\end{equation}
Substituting \eqref{614-1} and \eqref{615-1} into \eqref{614}, we obtain
\begin{equation}\label{615-2}
\begin{aligned}
	y_1^\prime=& \frac{\lambda \a\ln\phi}{2K_c}(1+t)^{-\lambda-1}
	-\frac{\a}{2(1+t)^\lambda}{A}\phi
	-K_d{A}^2\phi \\
	=&-K_cy_1^2-\frac{\a+2\a\ln\phi}{2(1+t)^\lambda}y_1
	-\frac{\a^2\ln\phi+\a^2(\ln\phi)^2}{4K_c(1+t)^{2\lambda}}
	+ \frac{\lambda \a\ln\phi}{2K_c}(1+t)^{-\lambda-1}\\
	=&-K_cy_1^2-\frac{\a+2\a\ln\phi}{2(1+t)^\lambda}y_1
	+\frac{2\lambda \a\ln\phi(1+t)^{\lambda-1}-\a^2\ln\phi(1+\ln\phi)}{4K_c(1+t)^{2\lambda}}\\
	:=&-a_2y_1^2-a_1y_1+a_0,
\end{aligned}
\end{equation}
where $a_0, a_1$ and $a_2$ are defined in \eqref{612}. This completes the proof of \eqref{610}.\\

For \eqref{611}, noticing that
$$
\begin{aligned}
	\phi^\backprime&=\phi_t-c\phi_x
	=\left(\sqrt{3K}\tau^{-1}\right)_t-c\phi_x\\
	&=-\sqrt{3K}\tau^{-2}\tau_t-c\phi_x=-cu_x-c\phi_x\\
	&=-c({w}-\phi)_x-c\phi_x=-c{w}_x+c\phi_x-c\phi_x\\
	&=-c{A},
\end{aligned}
$$
then using \eqref{65}, we have
$$
\begin{aligned}
	{B}^\backprime= -\frac{\a}{2(1+t)^\lambda}({B}-\frac{1}{c}\phi^\backprime)+K_d(-\frac{1}{c}\phi^\backprime{B}-{B}^2).
\end{aligned}
$$
Thus following a similar calculation in \eqref{614}-\eqref{615-2}, we obtain
%$$
%\begin{aligned}
	%q={B}\phi-\frac{\a}{2K_c(1+t)^\lambda}\ln\phi,
%\end{aligned}
%$$
%then
%$$
%\begin{aligned}
	%{B}=q\phi^{-1}+\frac{\a\ln\phi}{2K_c(1+t)^\lambda}\phi^{-1}.
%\end{aligned}
%$$
%We can get that $q$ satisfies
\begin{equation}\label{615-3}
\begin{aligned}
	q_1^\backprime=& \frac{\lambda \a\ln\phi}{2K_c}(1+t)^{-\lambda-1}
	-\frac{\a}{2(1+t)^\lambda}{B}\phi-K_d{B}^2\phi \\
	=&-K_cq_1^2-\frac{\a+2\a\ln\phi}{2(1+t)^\lambda}q_1
	-\frac{\a^2\ln\phi+\a^2(\ln\phi)^2}{4K_c(1+t)^{2\lambda}}
	+ \frac{\lambda \a\ln\phi}{2K_c}(1+t)^{-\lambda-1}\\
	=&-K_cq_1^2-\frac{\a+2\a\ln\phi}{2(1+t)^\lambda}q_1
	+\frac{2\lambda \a\ln\phi(1+t)^{\lambda-1}-\a^2\ln\phi(1+\ln\phi)}{4K_c(1+t)^{2\lambda}}\\
	:=&-a_2q_1^2-a_1q_1+a_0,
\end{aligned}
\end{equation}
where $a_2,$ $a_1,$ $a_0$ are defined in  (\ref{612}). Hence the proof of Lemma \ref{L6.2} is completed.\\

\textbf{Step 2.} For $\gamma=3$ and $\lambda\geq0$, we define
 $$
 h_1(\t) =\int_{\tilde{C}_0^{-1}}^{\tau}c(s)^{\frac{1}{2}}ds,
 $$
 where $\tilde{C}_{0}$ is the constant given in Lemma \ref{lem1.1}. It follows from \eqref{61-1} that $h_1(\tau)=(3K)^{\frac{1}{4}}\ln(\tilde{C}_0\tau)$.
 By replacing $h(\tau)$ with $h_1(\tau)$ and using the same calculations in the Steps 2-3 in subsection 3.1, we obtain (\ref{521}), that is,
$$
(3K)^{1\over 4}\tau^{-1}\tau_t+{\alpha\over {(1+t)^\lambda}}(3K)^{\frac{1}{4}}\ln(\tilde{C}_0\tau)\leq M,
$$
which implies that
%\begin{equation}\label{004.23}
%\tau^{-1}\tau_t+{\alpha\over {(1+t)^\lambda}}(\ln\tau-\ln\tilde{C}_0^{-1})\leq (3K)^{-\frac{1}{4}}M.
%\end{equation}
%By (\ref{004.23}), we have
$$
(\ln(\tilde{C}_{0}\tau))_t+{\alpha\over {(1+t)^\lambda}}\ln(\tilde{C}_{0}\tau)
\leq(3K)^{-\frac{1}{4}}M.
$$
Thus for $\lambda\neq 1$,
\begin{equation}\label{004.25}
	\left(e^{\frac{\alpha(1+t)^{1-\lambda}}{1-\lambda}}\ln(\tilde{C}_{0}\tau)\right)_t\leq  (3K)^{-\frac{1}{4}}Me^{\frac{\alpha(1+t)^{1-\lambda}}{1-\lambda}},
\end{equation}
and for $\lambda=1$,
\begin{equation}\label{004.25-1}
	\left((1+t)^{\alpha}\ln(\tilde{C}_{0}\tau)\right)_t\leq  (3K)^{-\frac{1}{4}}M(1+t)^{\alpha},
\end{equation}
where $M$ is the found given in \eqref{3.20-1}.

Integrating (\ref{004.25})-\eqref{004.25-1} from $0$ to $t$ respectively to see
$$
\begin{aligned}
e^{\frac{\alpha(1+t)^{1-\lambda}}{1-\lambda}}\ln(\tilde{C}_{0}\tau)
-e^{\frac{\alpha}{1-\lambda}}\ln(\tilde{C}_{0}\tau_0)
&\leq  (3K)^{-\frac{1}{4}}M\int_0^t e^{\frac{\alpha(1+s)^{1-\lambda}}{1-\lambda}}ds\\
&\leq  (3K)^{-\frac{1}{4}}Me^{\frac{\alpha(1+t)^{1-\lambda}}{1-\lambda}}t,\quad \quad \quad \text{for }\lambda\neq 1,
\end{aligned}
$$
and
$$
\begin{aligned}
	(1+t)^{\alpha}\ln(\tilde{C}_{0}\tau)
	-\ln(\tilde{C}_{0}\tau_0)
	&\leq  (3K)^{-\frac{1}{4}}M \int_0^t(1+s)^{\alpha}ds\\
	&\leq  (3K)^{-\frac{1}{4}}M(1+t)^{\alpha}t,\qquad\qquad\qquad \quad \text{for }\lambda=1.
\end{aligned}
$$
Therefore
$$
\begin{aligned}
\ln(\tilde{C}_{0}\tau)
&\leq \left[e^{\frac{\alpha}{1-\lambda}}\ln(\tilde{C}_{0}\tau_0)
+  (3K)^{-\frac{1}{4}}M e^{\frac{\alpha(1+t)^{1-\lambda}}{1-\lambda}}t\right]e^{-\frac{\alpha(1+t)^{1-\lambda}}{1-\lambda}}\\
&= e^{\frac{\alpha(1-(1+t)^{1-\lambda})}{1-\lambda}}\ln(\tilde{C}_{0}\tau_0)+  (3K)^{-\frac{1}{4}}Mt\\
&\leq\ln(\tilde{C}_{0}\tau_0)+(3K)^{-\frac{1}{4}}Mt,\quad \quad \quad \text{for }\lambda\neq 1,
\end{aligned}
$$
and
$$
\begin{aligned}
	\ln(\tilde{C}_{0}\tau)
	&\leq \left[\ln(\tilde{C}_{0}\tau_0)
	+  (3K)^{-\frac{1}{4}}M (1+t)^{\alpha}t\right](1+t)^{-\alpha}\\
	&= (1+t)^{-\alpha}\ln(\tilde{C}_{0}\tau_0)+  (3K)^{-\frac{1}{4}}Mt\\
	&\leq\ln(\tilde{C}_{0}\tau_0) +  (3K)^{-\frac{1}{4}}Mt,\qquad \quad \text{for }\lambda=1.
\end{aligned}
$$
Hence
$$
\begin{aligned}
\tau\leq\tau_0e^{(3K)^{-\frac{1}{4}}Mt},\quad \text{for }\lambda\geq 0,
\end{aligned}
$$
i.e.,
\begin{equation}\label{004.25-2}
\begin{aligned}
\rho\geq &\tau_0^{-1}e^{-(3K)^{-\frac{1}{4}}Mt}
\geq C_{0}^{-1}e^{-\hat{M}t},\quad \text{for }\lambda\geq 0,
\end{aligned}
\end{equation}
where we have used \eqref{1.10} and $\hat{M}:=(3K)^{-\frac{1}{4}}M$.\\

%%%%%%%%%%%%%%%%%%%%%%%%%%%%%%%%%%%%%%%%%%%%%%%%%%\gamma=3和\lambda=0的时候的密度下界的证明
%Especially, for $\gamma=3$ and $\lambda=0$, we obtain (\ref{521}), that is,
%$$
%(3K)^{1\over 4}\tau^{-1}\tau_t+\alpha(3K)^{\frac{1}{4}}\ln(\tilde{C}_0\tau)\leq M,
%$$
%which implies that
%$$
%(\ln(\tilde{C}_{0}\tau))_t+\alpha\ln(\tilde{C}_{0}\tau)
%\leq(3K)^{-\frac{1}{4}}M.
%$$
%Thus
%$$
%	\left(e^{\alpha(1+t)}\ln(\tilde{C}_{0}\tau)\right)_t\leq  (3K)^{-\frac{1}{4}}Me^{\alpha(1+t)},
%$$
%where $M$ is the found given in \eqref{3.20-1}.
%Integrating above inequality from $0$ to $t$ respectively to see
%$$
%\begin{aligned}
%e^{\alpha(1+t)}\ln(\tilde{C}_{0}\tau)
%-e^{\alpha}\ln(\tilde{C}_{0}\tau_0)
%&\leq  (3K)^{-\frac{1}{4}}M\int_0^t e^{\alpha(1+s)}ds\\
%&= \frac{(3K)^{-\frac{1}{4}}M}{\alpha}(e^{\alpha(1+t)}-e^{\alpha})\\
%&\leq  \frac{(3K)^{-\frac{1}{4}}M}{\alpha}e^{\alpha(1+t)},
%\end{aligned}
%$$
%Therefore
%$$
%\begin{aligned}
%\ln(\tilde{C}_{0}\tau)
%&\leq \left[e^\alpha\ln(\tilde{C}_{0}\tau_0)
%+ \frac{(3K)^{-\frac{1}{4}}M}{\alpha}e^{\alpha(1+t)}\right]e^{-{\alpha(1+t)}}\\
%&= e^{\alpha(1-(1+t))}\ln(\tilde{C}_{0}\tau_0)+\frac{(3K)^{-\frac{1}{4}}M}{\alpha}\\
%&\leq\ln(\tilde{C}_{0}\tau_0)+\frac{(3K)^{-\frac{1}{4}}M}{\alpha},
%\end{aligned}
%$$
%Hence
%$$
%\begin{aligned}
%\tau\leq\tau_0e^{\frac{(3K)^{-\frac{1}{4}}M}{\alpha}},
%\end{aligned}
%$$
%i.e.,
%$$
%\begin{aligned}
%\rho\geq &\tau_0^{-1}e^{-\frac{(3K)^{-\frac{1}{4}}M}{\alpha}}
%\geq C_{0}^{-1}e^{-\hat{M}},
%\end{aligned}
%$$
%where we have used \eqref{1.10} and $\hat{M}:=\frac{(3K)^{-\frac{1}{4}}M}{\alpha}$.\\
%%%%%%%%%%%%%%%%%%%%%%%%%%%%%%%%%%%%%%%

\textbf{Step 3.} For $\gamma=3$, $\lambda<0$, By replacing $h(\tau)$ with $h_1(\tau)$ and using the same calculations in the Step 4 in subsection 3.1, we can obtain (\ref{520}), that is,
$$
(3K)^{1\over 4}\tau^{-1}\tau_t+{\alpha\over {(1+t)^\lambda}}(3K)^{\frac{1}{4}}\ln(\tilde{C}_0\tau)\leq e^{-\frac{\lambda(1+t)^{1-\lambda}}{1-\lambda}}M,
$$
which implies that
%\begin{equation}\label{0004.23}
%\tau^{-1}\tau_t+{\alpha\over {(1+t)^\lambda}}\big(\ln\tau-\ln\tilde{C}_0^{-1}\big)\leq (3K)^{-\frac{1}{4}}e^{-\frac{\lambda(1+t)^{1-\lambda}}{1-\lambda}}M.
%\end{equation}
%By (\ref{0004.23}), we have
$$
(\ln(\tilde{C}_{0}\tau))_t+{\alpha\over {(1+t)^\lambda}}\ln(\tilde{C}_{0}\tau)
\leq(3K)^{-\frac{1}{4}}e^{-\frac{\lambda(1+t)^{1-\lambda}}{1-\lambda}}M.
$$
Thus
\begin{equation}\label{0004.25}
\begin{aligned}
\left(e^{\frac{\alpha(1+t)^{1-\lambda}}{1-\lambda}}\ln(\tilde{C}_{0}\tau)\right)_t
\leq\left((3K)^{-\frac{1}{4}}e^{-\frac{\lambda(1+t)^{1-\lambda}}{1-\lambda}}M\right)e^{\frac{\alpha(1+t)^{1-\lambda}}{1-\lambda}}.
\end{aligned}
\end{equation}
Integrating (\ref{0004.25}) from $0$ to $t$ to see
$$
\begin{aligned}
&\ \ \ \ e^{\frac{\alpha(1+t)^{1-\lambda}}{1-\lambda}}\ln(\tilde{C}_{0}\tau)-e^{\frac{\alpha}{1-\lambda}}\ln(\tilde{C}_{0}\tau_0)\\
&\leq (3K)^{-\frac{1}{4}}M\int_0^te^{\frac{(\alpha-\lambda)(1+s)^{1-\lambda}}{1-\lambda}}ds\\
&= (3K)^{-\frac{1}{4}}Me^{\frac{(\alpha-\lambda)(1+t)^{1-\lambda}}{1-\lambda}}t.
\end{aligned}
$$
Therefore
$$
\begin{aligned}
\ln(\tilde{C}_{0}\tau)
&\leq \left[e^{\frac{\alpha}{1-\lambda}}\ln(\tilde{C}_{0}\tau_0)
+ (3K)^{-\frac{1}{4}}Me^{\frac{(\alpha-\lambda)(1+t)^{1-\lambda}}{1-\lambda}}t\right]e^{-\frac{\alpha(1+t)^{1-\lambda}}{1-\lambda}}\\
&= e^{\frac{\alpha(1-(1+t)^{1-\lambda})}{1-\lambda}}\ln(\tilde{C}_0\tau_0) +  (3K)^{-\frac{1}{4}}Me^{\frac{-\lambda(1+t)^{1-\lambda}}{1-\lambda}}t\\&
\leq\ln\tilde{C}_0\tau_0+(3K)^{-\frac{1}{4}}Me^{\frac{-\lambda(1+t)^{1-\lambda}}{1-\lambda}}t\\
&=\ln(\tilde{C}_{0}\tau_0)+(3K)^{-\frac{1}{4}}Me^{\frac{-\lambda(1+t)^{1-\lambda}}{1-\lambda}}t.
\end{aligned}
$$
Hence
$$
\begin{aligned}
\tau\leq\tau_0e^{(3K)^{-\frac{1}{4}}M t e^{\frac{-\lambda(1+t)^{1-\lambda}}{1-\lambda}}},\quad \text{for }\lambda<0,
\end{aligned}
$$
i.e.,
\begin{equation}\label{0004.25-1}
\begin{aligned}
\rho \geq \tau_0^{-1}e^{-(3K)^{-\frac{1}{4}}M t e^{\frac{-\lambda(1+t)^{1-\lambda}}{1-\lambda}}}
\geq C_{0}^{-1}e^{-\hat{M} t e^{\frac{-\lambda(1+t)^{1-\lambda}}{1-\lambda}}},\quad \text{for }\lambda<0,
\end{aligned}
\end{equation}
where we have used \eqref{1.10} and $\hat{M}:=(3K)^{-\frac{1}{4}}M$..\\
%\begin{equation}\label{0e.3}
%M_3=\tau_0^{-1},\ \  \tilde{M}_4=(3K)^{-\frac{1}{4}}M.
%\end{equation}

%\begin{lemma}\label{L6.1}
	%The smooth solutions  of (\ref{1.1})-(\ref{1.2}) satisfy
	%\begin{equation}\label{65}
		%\begin{aligned}
			%{A}^\prime= -\frac{\a}{2(1+t)^\lambda}({A}+{B})+K_d({A}{B}-{A}^2)
		%\end{aligned}
%\end{equation}
%and
	%\begin{equation}\label{66}
		%\begin{aligned}
			%{B}^\backprime=-\frac{\a}{2(1+t)^\lambda}({A}+{B})+K_d({A}{B}-{B}%^2),
		%\end{aligned}
	%\end{equation}
	%where $K_d = K_c\phi$,
	%${A}={w}_x,$ ${B}={z}_x$.
%\end{lemma}

%\textbf{Proof}. Take the partial derivative of both sides of equation (\ref{63}) with respect to $x$,
%\begin{equation}\label{67}
%\begin{aligned}
%({w}^\prime)_x=-\frac{\a}{2(1+t)^\lambda}({z}+{w})_x,
 %\end{aligned}
%\end{equation}
%then
%\begin{equation*}\label{20.11}
%\begin{aligned}
%({w}^\prime)_x=({w}_t+c{w}_x)_x={w}_{tx}+cw_{xx}+c_{x}{w}_{x}=({w}_x)^\prime+c_{x}{w}_{x}
%=-\frac{\a}{2(1+t)^\lambda}({z}_x+{w}_x).
 %\end{aligned}
%\end{equation*}
%Thus,
%\begin{equation*}\label{20.12}
%\begin{aligned}
%{A}^\prime+c_{x}{A}=-\frac{\a}{2(1+t)^\lambda}({A}+{B}).
% \end{aligned}
%\end{equation*}
%And
%\begin{equation*}\label{20.13}
%\begin{aligned}
%c_x = \left(K_c\phi^{2}\right)_x
%=2 K_c\phi \phi_x=2K_d\phi_x=K_d({A}-{B}),
% \end{aligned}
%\end{equation*}
%where we used (\ref{62}), and denote $K_d=K_c\phi$.
%Thus (\ref{67}) changes to
%\begin{equation*}\label{20.15}
%\begin{aligned}
%{A}^\prime+K_d({A}-{B}){A}= -\frac{\a}{2(1+t)^\lambda}({A}+{B}),
 %\end{aligned}
%\end{equation*}
%and we get (\ref{65}). The calculation of (\ref{66}) is same.

\textbf{Step 4.} For any given $t_0>0$, it follows from \eqref{61}, \eqref{612}, \eqref{004.25-2} and \eqref{0004.25-1} that
\begin{equation}\label{9.9}
\begin{aligned}
	&|a_1|\leq \bigg{[}\frac{\alpha}{2}+\alpha\ln \bar{M}(t_0)\bigg{]}\max\{(1+t_0)^{-\lambda},1\},\\
&|a_0|\leq \bigg{[}\frac{|\lambda|\alpha}{2K_{c}}\ln \bar{M}(t_0)+\frac{\alpha^2}{4K_{c}}\ln \bar{M}(t_0)+\frac{\alpha^2}{4K_{c}}(\ln \bar{M}(t_0))^2\bigg{]}\max\{(1+t_0)^{-2\lambda},1\},
\end{aligned}\end{equation}
where for $\lambda\geq 0$,
\begin{equation}\label{9.9-2}
\bar{M}(t_0):=\max\big\{\sqrt{3K}\tilde{C}_{0}, (\sqrt{3K}\tilde{C}_{0})^{-1}, \sqrt{3K}e^{-\hat{M}t_0}, (\sqrt{3K}e^{-\hat{M}t_0})^{-1}\big\}\geq 1,
\end{equation}
and for $\lambda<0$,
\begin{equation}\label{9.9-3}
\bar{M}(t_0):=\max\big\{\sqrt{3K}\tilde{C}_{0}, (\sqrt{3K}\tilde{C}_{0})^{-1}, \sqrt{3K}e^{-\hat{M}t_0e^{\frac{-\lambda}{1-\lambda}(1+t_0)^{-\lambda}}}, (\sqrt{3K}e^{-\hat{M}t_0e^{\frac{-\lambda}{1-\lambda}(1+t_0)^{-\lambda}}})^{-1}\big\}\geq 1.
\end{equation}
It is directly to check that when $\bar{M}(t_0)$ is non-decreasing with respect to $t_0$.

Using \eqref{612}, \eqref{9.9} and \eqref{9.9-2}-\eqref{9.9-3}, there exists a function $\tilde{M}(t_0)$, which is non-decreasing with respect to $t_0$, such that
\begin{equation}\label{9.9-1}
\left\vert\frac{a_1}{a_2}\right\vert+\sqrt{\left\vert\frac{a_1}{a_2}\right\vert^2+2\left\vert\frac{a_0}{a_2}\right\vert}\leq \tilde{M}(t_0),\quad \text{for any }0\leq t\leq t_0.
\end{equation}
If there exists one point $x_0$ such that $y_1(x_0, 0)<-\max\big\{\frac{2}{K_{c}t_0}, \tilde{M}(t_0)\big\}$, i.e., \eqref{9.14} holds, then using \eqref{9.9-1}, one has
$$
-\frac{1}{2}a_2y_1(x_0, 0)^2-a_1(0)y_1(x_0, 0)+a_0(0)< 0,
$$
and using \eqref{610}, we obtain
$$
y_1^{\prime}\vert_{t=0}=-\frac{1}{2}a_2y_1(x_0, 0)^2+\left(-\frac{1}{2}a_2y_1(x_0, 0)^2-a_1(0)y_1(x_0, 0)+a_0(0)\right)< -\frac{1}{2}a_2y_1(x_{0},0)^2\leq 0.
$$
Thus along the forward characteristic line starting from $(x_0,0)$, we have
\begin{equation}\label{9.10}
y_1(x_{+}(t),t)<y_1(x_0,0)< -\max\big\{\frac{2}{K_{c}t_0}, \tilde{M}(t_0)\big\}, \quad \text{for all }0\leq t\leq t_0.
\end{equation}
Turning back to \eqref{610} and using \eqref{9.9-1} and \eqref{9.10}, it holds
\begin{equation}\label{9.11}
y_1^{\prime}=-\frac{1}{2}a_2y_1^2+\left(-\frac{1}{2}a_2y^2-a_1y+a_0\right)<-\frac{1}{2}a_2y_1^2,\quad \text{for all }0\leq t\leq t_0.
\end{equation}
Integrating \eqref{9.11} with respect to $t$ along the forward characteristic line starting from $(x_0,0)$, we get
\begin{equation}\label{9.12}
	0>y_1^{-1}(x_{+}(t),t)> y_1^{-1}(x_0,0)+\int_{0}^{t}\frac{1}{2}a_2ds=y_1^{-1}(x_0,0)+\frac{1}{2}K_{c}t,\quad \text{for all }0\leq t\leq t_0.
\end{equation}
Since $y_1^{-1}(x_0,0)\in (-\frac{1}{2}K_{c}t_0,0)$, the right hand side of \eqref{9.12} will goes to zero as $t\rightarrow t_0$, which implies that $y(x_{+}(t),t)$ will blow up before $t_0$. By similar arguments, we obtain that $q_1$ will blow up before $t_0$ if \eqref{9.15} holds.

Finally, using \eqref{004.25-2}, \eqref{0004.25-1} and \eqref{611}, we see that $y_1$ or $q_1$ blows up in finite time means $\|(\tau_{x},u_{x})\|_{L^{\infty}}$ blows up in finite time. Therefore the proof of Theorem \ref{thm9.1} is completed.
$\hfill\square$

\end{document}